\documentclass[12pt]{amsart}  
\usepackage{amsmath, amscd, amssymb}

\newcommand{\mybf}{}
\newcommand{\myem}{}

\newtheorem{thm}{Theorem}[section]
\newtheorem{lem}[thm]{Lemma}

\newtheorem{remark}[thm]{Remark}
\newtheorem{cor}[thm]{Corollary}
\newtheorem{prop}[thm]{Proposition}

\theoremstyle{definition}
\newtheorem{example}[thm]{Example}
\newtheorem{construction}[thm]{Construction}
\newtheorem{defn}[thm]{Definition}

\newcommand{\M}{\mathcal{M}}
\newcommand{\CO}{\mathcal{O}}

\newcommand{\U}{\mathcal{U}}

\newcommand{\Grass}{\mathrm{Grass}}
\newcommand{\Map}{\mathrm{Map}}
\newcommand{\THH}{\mathrm{THH}}
\newcommand{\os}{\mathcal{I}_S}
\newcommand{\osp}{\mathrm{Sp}\!^O}

\newcommand{\ev}{\mathrm{ev}}
\newcommand{\ho}{\mathrm{ho}}

\newcommand{\So}{\mathbb{S}}

\newcommand{\F}{\mathcal{F}}
\newcommand{\FG}{\mathcal{F}_G}

\newcommand{\supp}{\mathrm{supp}}

\newcommand{\map}{{\mathrm{map}}}
\newcommand{\mor}{{\mathrm{mor}}}

\newcommand{\Cat}{{\rm{Cat}}}
\newcommand{\Fun}{{\rm{Fun}}}
\newcommand{\GR}{{\mathbb{G}}}
\newcommand{\Gr}{{\mathcal{G}}}
\newcommand{\W}{{\mathbb{W}}}

\newcommand{\Z}{{\mathbb{Z}}}
\newcommand{\CC}{{\mathbb{C}}}
\newcommand{\C}{{\mathcal{C}}}
\newcommand{\D}{{\mathcal{D}}}
\newcommand{\N}{{\mathbb{N}}}
\newcommand{\op}{{\operatorname{op}}}
\newcommand{\pr}{\operatorname{pr}}
\newcommand{\ob}{\operatorname{ob}}
\newcommand{\id}{\operatorname{id}}
\newcommand{\sk}{\mathrm{sk}}
\newcommand{\Topp}{{\mathcal{T}}}

\newcommand{\Set}{\mathcal{E}\! \mathit{ns}}
\newcommand{\Ab}{\mathcal{A}\! \mathit{b}}

\newcommand{\Ch}{\mathcal{C}\! \mathit{h}}
\newcommand{\Ens}{\mathcal{F}\!\mathit{in}}
\newcommand{\colim}[1]{\mathop{\underset{#1}
            {{\text{\rm colim}}}}}

\author{M. Brun}
\title{Witt Vectors and Equivariant Ring Spectra}
\begin{document}
\begin{abstract}
This paper establishes a connection between equivariant ring spectra
and Witt vectors in the sense of Dress and Siebeneicher. Given a
commutative ringspectrum $T$ in the highly structured sense, that is, an
$E_{\infty}$-ringspectrum, with action of a finite group $G$ we construct a
ringhomomorphism from the ring of $G$-typical Witt vectors of the
zeroth homotopy group of 
$T$ to the zeroth homotopy group of the $G$-fixed point spectrum of
$T$. In the particular case, where $T$ is the periodic unitary
cobordism spectrum 
introduced by Strickland, we show that this ringhomomorphism is
injective, and we 
interpret this in terms of equivariant cobordism.   
\end{abstract}
\maketitle
\section{Introduction}
Tambara has studied the interaction of trace,
norm and restriction maps connecting Burnside rings of different subgroups
of a finite group $G$. He observed, that certain relations
between the 
trace-, norm- and restriction maps satisfied by these Burnside
rings are also satisfied by the trace-, norm- and restriction maps 
between even cohomology groups- and between representation rings of
subgroups of $G$. A
collection of abelian groups with trace-, norm- and
restriction maps satisfying such
relations was called a $\mathrm {TNR}$-functor in \cite{Tambara}.
In \cite{Brun} we called such a collection a $G$-Tambara functor.
Here we call it
simply a Tambara functor. 

The main result of this paper states that
every $E_\infty$ ring spectrum $T$ with an action of $G$ has an associated
Tambara functor $\widetilde T$ with
$$\widetilde T(X) = [\Sigma^{\infty}X_+ ,T]_G.$$
Here $\Sigma^{\infty} X_+$ denotes the suspension spectrum on $X$ with
a disjoint added basepoint and $[\Sigma^{\infty}X_+ ,T]_G$ denotes the
morphisms from $\Sigma^{\infty}X_+$ to $T$ in the equivariant stable
homotopy category associated to a complete universe of
$G$-representations in the sense of Lewis--May--Steinberger \cite{LMS}.

The main result of \cite{Brun} is that 
there is a homomorphism $\tau_S \colon \W_G(S(G/e)) \to
S(G/G)$ 
for every Tambara
functor $S$. 
Here $\W_G(S(G/e))$ is Dress and Siebeneicher's
ring of $G$-typical Witt vectors 
on the commutative ring $S(G/e)$ introduced in \cite{Dress-Siebeneicher}. 
Together these two results give a homomorphism
$$\tau_{\widetilde T} \colon \W_G([\So,T]) \to [\So,T]_G$$
for every $E_\infty$ ring
spectrum $T$ with an action of 
$G$. 
Here $\So$ denotes the 
sphere spectrum and we have
identified $[\Sigma^{\infty}G_+,T]_G$ with $[\So,T]$ and we have
identified $\Sigma^{\infty}(G/G)_+$ with $\So$. The group $[\So,T]$ is
the zeroth homotopy group $\pi_0(T)$ of $T$ and the group $[\So,T]_G$ is the
zeroth homotopy group $\pi_0(T^G)$ of the $G$-fixed point spectrum $T^G$
associated to $T$.
The following theorem of Dress and Siebeneicher
\cite{Dress-Siebeneicher} 
shows that in
 some cases $\tau_{\widetilde T}$ gives
information on
$[\So,T]_G$. 
\begin{thm}
  The homomorphism 
  $\tau_{\widetilde \So}
  : \W_G([\So,\So]) \to [\So,\So]_G$ 
  is an isomorphism.
\end{thm}
Hesselholt and Madsen have extended this result  
to topological Hochschild homology
in the special case where $G$ is a finite cyclic group
\cite{Hesselholt-Madsen}. Given an $E_{\infty}$ ring-spectrum $R$,
$\THH(R)$ denotes the topological Hochschild homology of $R$
considered as an $E_\infty$ ring-spectrum with action of the circle group.
\begin{thm}
  If $G$ a finite cyclic group and $R$ is an $E_\infty$
  ring spectrum,
  then the homomorphism
  $\tau_{\widetilde{\THH(R)}} : \W_G([\So,\THH(R)]) \to
  [\So,\THH(R)]_G$ is an isomorphism. 
\end{thm}
Let $MP$ denote the $G$-equivariant $E_\infty$ ring
spectrum representing the periodic version of the unitary Thom spectrum  
introduced by Strickland \cite{Strickland}. 
(This spectrum is intimately related to cobordism of $G$-manifolds.)
For $H \leq G$ the group $[\So,MP]_H$ is the sum of
the even dimensional homotopy groups of $MU^H$. Using that the
spectrum $MP$ has 
a restriction map similar to the restriction map for topological
Hochschild homology we prove: 
\begin{thm}
\label{indthmfirst}
  The homomorphism 
  $\tau_{\widetilde{MP}}
  : \W_G([\So,MP]) \to [\So,MP]_G$
  is injective.
\end{thm}
The homomorphism $\tau_{\widetilde{MP}}$ is not an isomorphism. Already
in the case where $G$ is a cyclic group of prime 
order it is not surjective.

The equivariant unitary cobordism ring $\U^G_*$ of $G$-manifolds
is related to the non-equivariant cobordism ring $\U_{*}
\cong [\So,MP]$ as follows: 
\begin{thm}
\label{indthmlast}
  The ring $\W_G(\U_{*})$ embeds as a subring of
  $\U^G_{*}$.
\end{thm}
In order to construct the Tambara functor $\widetilde T$ associated to
an $E_\infty$-ring spectrum $T$ with action of $G$ we need to study ordinary
induction  
and 
smash-induction of equivariant spectra. 
Since ordinary induction is constructed by
wedge sums and smash-induction is constructed by smash products,
the interaction
between smash-products and wedge-sums will play an important role for
us. One way to encode the relations between wedge sums and
smash products is in Laplaza's concept of a bimonoidal category
\cite{Laplaza} defined by a large number of commutative diagrams. 
In order to apply Laplaza's coherence result we shall reformulate it
in such a way that it becomes clear that commutative ring-objects in
bimonoidal categories give rise to Tambara functors. 
This reformulation involves a new approach to higher coherences
via partial pseudo-functors. 

For the passage from categorical constructions to homotopy groups we
need a homotopical analysis of the 
smash induction of spectra. 
In particular we must control the
operation taking a (non-equivariant) spectrum $X$ to its $G$'th smash-power
$X^{\wedge G}$ with $G$ acting by permuting the smash factors. 
This is gained by examining the interaction between
cofibrantions and smash induction.   
Smash induction is sensitive to the choice of
category of equivariant spectra. In this paper we have chosen to  
work with the orthogonal
spectra of \cite{MMSS} and \cite{MM}. 
We warn the reader that the smash-induction of
an $\So$-modules in the sense of \cite{EKMM} 
might have equivariant 
homotopy type different from the smash-induction of its associated
orthogonal spectrum.
However the results of this paper also apply to the category of
equivariant symmetric spectra with the
model structure considered in \cite{HSS}, \cite{H-symm} and
\cite{M-symm}. 
In a joint paper in preparation with M. Lydakis we show that they also
apply to 
the category of equivariant Gamma-spaces. 

For our construction of the Tambara functor
$\widetilde T$
we actually need $T$ to be 
a strictly commutative orthogonal ring spectrum with action of
$G$. However it is well-known that every orthogonal spectrum with an
action of an $E_\infty$-operad is stably equivalent to a strictly
commutative ring 
spectrum.
(See
\cite[Lemma III.8.4]{MM}, \cite[Remark 0.14]{MMSS} and the proof of
\cite[Proposition II.4.3]{EKMM}). 
We cicumvent technical difficulties caused by the fact that in general the
commutative replacement of $T$ can not be chosen to be 
cofibrant in the model structure considered in this paper.

The paper is organized as follows:
In Section \ref{tamfunctors} we give some preliminaries on Tambara
functors. Section \ref{tamcatsec} introduces 
Tambara categories. These are category-valued coherent Tambara
functors playing a key role for the 
construction of Tambara functors from commutative ring objects in
bimonoidal categories.
In section \ref{indofcof} we give a criterion for cofibrations to be
preserved under smash induction. This section is the main technical
part of the construction of $\widetilde T$ from $T$.
Section \ref{firsttrfsec} explains how to construct  Tambara functors
from commutative ring objects in a bimonoidal category.
In Section \ref{sechtytam} we show that the category of orthogonal
spectra satisfies the criteria of Section \ref{indofcof}. As an aside
we show that there is a Tambara category of chain complexes.  
Section \ref{sechotrf} ends the construction of 
$\widetilde T$.
In Section \ref{witam} we study the equivariant version of
Strickland's spectrum $MP$ and prove Theorem
\ref{indthmfirst} and Theorem \ref{indthmlast}.   
Section \ref{secfilob} is an appendix on filtered
objects in symmetric monoidal categories.

\section{Tambara Functors}
\label{tamfunctors}
For the convenience of the reader we recollect parts of the work
\cite{Tambara} of 
D. Tambara.
The category $\FG$ of 
finite left $G$-sets 
is small, contains finite sums, finite products and
  pullbacks.
In particular $\FG$ contains an initial object $\emptyset$
  and a final object 
  $*$.
It is well-known that 
  for every  $f : X \to Y$ in $\FG$, the pull-back functor
  \begin{eqnarray*}
    \FG/Y &\to& \FG/X, \\
    (B \to Y) &\mapsto& (X \times_Y B \to X)
  \end{eqnarray*}  
  has a right adjoint
  \begin{eqnarray*}
    \Pi_f : \FG/X &\to& \FG/Y \\
    (A \xrightarrow p X) &\mapsto& (\Pi_fA \xrightarrow{\Pi_f p} Y).
  \end{eqnarray*}
The $G$-map $\Pi_f p$ has fibers $(\Pi_f p)^{-1}(y) = \prod_{x \in
f^{-1}(y)} p^{-1}(x)$. 
\begin{defn}
\label{defexpdiag}
  A diagram in $\FG$ isomorphic to a diagram of the form:
\begin{displaymath}
\begin{CD}
  X @<p<< A @<{e}<< X \times_Y \Pi_fA \\
  @VfVV @. @VV{f'}V \\
  Y @= Y @<{{\Pi_f p}}<< \Pi_f A,
\end{CD}
\end{displaymath}
where $f'$ is the projection and $e$ is adjoint
to the identity on 
$\Pi_f A$  is called an {\em exponential diagram}.
\end{defn}
Two diagrams $X \leftarrow A \to B \to Y$ and $X \leftarrow
A' \to B' \to Y$ in $\FG$ are {\em equivalent} if there exist
isomorphisms $A \to A'$, $B \to B'$ in $\FG$ making the diagram
\begin{displaymath}
\begin{CD}
  X @<<< A @>>> B @>>> Y \\
  @| @VVV @VVV @| \\
  X @<<< A' @>>> B' @>>> Y   
\end{CD}
\end{displaymath}
commutative.
The category $U^G_+$ has finite $G$-sets as
objects, morphisms from $X$ to $Y$ given by equivalence classes
$[X \leftarrow A \to B \to Y]$ of diagrams $X \leftarrow A \to B \to
Y$  
and composition $\circ : U^{G}_+(Y,Z) \times U^{G}_+(X,Y)
\to U^{G}_+(X,Z)$ defined by the formula
\begin{displaymath}
  [Y \leftarrow C \to D \to Z] \circ [X \leftarrow A \to B \to Y] = [X
  \leftarrow A'' \to \widetilde D \to Z],  
\end{displaymath}
where the maps on the right are composites of the maps in the
diagram 
\begin{displaymath}
\begin{CD}
  X @<<< A @<<< A' @<<< A'' \\
  @. @VVV @VVV @VVV \\
  @. B @<<< B' @<<< \widetilde C \\
  @. @VVV @VVV @| \\
  @. Y @<<< C @. \widetilde C \\
  @. @. @VVV @VVV \\
  @. Z @<<< D @<<< \widetilde D.
\end{CD}
\end{displaymath}
Here the three squares are pull-back diagrams and the diagram
\begin{displaymath}
\begin{CD}
  C @<<< B' @<<< \widetilde C \\
  @VVV @. @VVV \\
  D @<<< \widetilde D @= \widetilde D
\end{CD}
\end{displaymath}
is an exponential diagram.
Given $f : X \to Y$ in $\FG$ we denote be $R_f,T_f$ and $N_f$ the morphisms
\begin{eqnarray*}
  R_f &=& [Y \xleftarrow f X \xrightarrow = X \xrightarrow = X], \\
  T_f &=& [X \xleftarrow = X \xrightarrow = X \xrightarrow f Y], \\
  N_f &=& [X \xleftarrow = X \xrightarrow f Y \xrightarrow = Y]
\end{eqnarray*}
in $U^G_+$.
Every morphism in $U^{G}_+$ can be written as a composition of
morphisms on the form $R_f$, $T_f$ and $N_f$.
\begin{prop}\label{semi-ring-str}
(i) For every sum diagram $X_1 \xrightarrow {i_1} X \xleftarrow {i_2}
  X_2$ in $\FG$, the diagram
$X_1 \xleftarrow {R_{i_1}} X \xrightarrow {R_{i_2}} X_2$ is a product
diagram in $U^{G}_+$. The object $\emptyset$ is final in $U^{G}_+$.

(ii) Let $\nabla : X \amalg X \to X$ denote the fold morphism of an object $X$
of $\FG$ and let $i : \emptyset \to X$. Then $X$, considered as an object of
$U^G_+$,
is a semi-ring object with addition $T_{\nabla}$, additive unit $T_i$,
multiplication $N_{\nabla}$ and multiplicative unit $N_i$.

(iii) If $f : X \to Y$ is a morphism in $\FG$, then the morphisms
$R_f$, $T_f$ and 
$N_f$ of $U^{G}_+$ preserve the above structures of commutative
semi-ring, additive 
monoid and 
multiplicative monoid on $X$ and $Y$, respectively.
\end{prop}
\begin{defn}\label{tamdef}
  The category of {\em semi Tambara functors} is the category of
  set-valued product-preserving functors on $U^G_+$ with morphisms
  given by natural transformations.
\end{defn}
\begin{defn}\label{tamfundef}
  The category of {\em Tambara functors} is the full subcategory of the
  category of semi Tambara functors consisting of those
  product-preserving functors $S : U^G_+ \to \Set$ satisfying that the
  underlying additive monoid of $S(X)$ is an abelian group for
  every finite $G$-set $X$.
\end{defn}
Given a semi Tambara functor $S$ and $\phi \in
U^{G}_+(X,Y)$ we obtain a function $S(\phi) : S(X) \to S(Y)$. Since
$S$ is product-preserving, 
it follows from (ii) of Proposition \ref{semi-ring-str} that $S(X)$ is a
semi-ring. 
Given a morphism $f : X \to Y$ in $\FG$
we let $S^*(f) = S(R_f)$, $S_+(f) = S(T_f)$ and
$S_\bullet(f) = S(N_f)$.
It follows from (iii) of \ref{semi-ring-str} that $S^*(f)$ is a
homomorphism of semi-rings, that 
$S_+(f)$ is an additive homomorphism and that $S_\bullet(f)$ is
multiplicative. A semi Tamara functor $S$ is uniquely determined by the
functions $S^*(f)$, $S_+(f)$ and $S_\bullet(f)$ for $f$ in
$\FG$. 
\begin{example}
  A commutative ring with an action of $G$ gives rise to a Tambara
functor with value $\Map_G(X,R)$ on a finite $G$-set $X$, as explained
below in Example \ref{invariant-ring} and in \cite[Example 3.1]{Tambara}. 
\end{example}
\begin{example}
  Given a Tambara functor $S$, the commutative ring $S(G/e)$
  has an obvious 
  action of $G$. The forgetful functor from the category of Tambara
  functors to the category of commutative rings with action of
  $G$ has a left adjoint functor $L$. The commutative ring $L(R)(G/G)$
  is 
  closely related to the generalized Witt-vectors of Dress and
  Siebeneicher \cite{Dress-Siebeneicher}. In fact it was shown in
  \cite{Brun}, that if $G$ acts 
  trivially on $R$, then $L(R)(G/G)$ is isomorphic to the ring
  $\W_G(R)$ of $G$-Witt vectors on $R$. 
\end{example}
\begin{example}
  Let $A$ be an $E_\infty$-ring spectrum with an action of $G$. In
  Section \ref{sechotrf} we show that there is a Tambara functor 
  $\widetilde A$
  with $\widetilde A(X) = [\Sigma^{\infty} X_+,A]_G$. The special case
  where $A$ is the 
  sphere spectrum is \cite[Example 3.2]{Tambara}. In Section \ref{witam} we
  examine the case $A = MP$.
\end{example}
\section{Tambara Categories}
\label{tamcatsec}
Many Tambara functors are induced by lax
category-valued Tambara functors. In order to explain what a lax
category-valued Tambara functor is, it will be helpful to consider
partial categories as defined below and to consider a particular partial
category 
$sU_+$ codifying Laplaza's coherence result for bimonoidal
categories \cite{Laplaza}. 
Some authors call such categories bimonoidal categories.
\begin{defn}
  A {\em partial category} $\Gr$ consists of a class $\ob \Gr$ of objects, a
  class $\mor \Gr$ of arrows, domain and codomain functions $d,c : \mor
  \Gr \to \ob \Gr$, an identity function $\id : \ob \Gr \to \mor \Gr$, a
  subclass $\mathrm{com} (\Gr) \subseteq \mor \Gr \times_{\ob \Gr} \mor \Gr$
  of {\em composable arrows} and a composition $\circ : \mathrm{com} (\Gr) \to
  \mor \Gr$ subject to the following associativity and unit axioms:
  \begin{description}
  \item[Associativity] For a diagram
        $  a \xrightarrow f b \xrightarrow g c \xrightarrow k d$
        of composable arrows in $\Gr$
        the relation
	$  k \circ (g \circ f) = (k \circ g) \circ f$ holds,
  \item[Unit] For arrows $f : a \to b$ and $g : b \to c$
        the pairs $(f,\id_b)$ and $(\id_b,g)$ of arrows are
        composable with 
	$\id_b \circ f = f$ and $g \circ \id_b = g.$    
  \end{description}
\end{defn}
\begin{defn}
  A {\em partial functor} $F : \C \to \D$ of partial
  categories consists of functions $F : \ob \C \to \ob \D$ and
  $F : \mor \C \to \mor \D$ compatible with the domain and codomain
  functions for $\C$ and $\D$ taking every pair $(f,g)$ of composable arrows
  of $\C$ to a pair $(F(f),F(g))$ 
  of composable arrows of $\D$ with 
  $F(g) \circ F(f) = F(g \circ f)$. 
\end{defn}
\begin{defn}
  A natural transformation $t : F \to G$ between partial
  functors $F,G : \C \to \D$ consists of arrows $t_c : F(c) \to G(c)$
  in $\D$ for every object $c$ of $\C$ subject to the condition that
  for every arrow $f : a \to b$ of $\C$ the pairs 
  $(G(f),t_a)$ and $(t_b,F(f))$ of arrows in $\D$ are composable with
  $G(f) \circ t_c = 
  t_b \circ F(f)$.
\end{defn}
\begin{defn}
An object $c$ of a partial category $\C$ is the product of the objects
$a$ and $b$ if there exist arrows $p_a: c \to a$ and $p_b: c \to b$
with the property
that given arrows $f: x \to a$ and $g : x \to b$ there exists a
unique arrow $h : x \to c$ such that the pairs of arrows $(p_a,h)$
and $(p_b,h)$ are composable with $p_a \circ h = f$ and $p_b
\circ h = g$. 
\end{defn}
\begin{remark}
  Products in partial categories are unique up to a unique
  isomorphism, and the notions of general limits and colimits make
  perfect sense in partial categories.
\end{remark}
\begin{remark}
  To a partial category $\Gr$ we can associate a directed graph
  $\widetilde \Gr$ with
  one arrow $(g,f)$ for each pair $a \xrightarrow f b \xrightarrow g
  c$ of composable
  morphisms in $\Gr$ together with the equivalence relation $\sim$ on the set
  of arrows generated by $(g,f)
  \sim (g \circ f, 
  \id_a) \sim (\id_b,g \circ f)$.
  There is an associated partial functor $i_{\Gr}: \Gr \to
  \widehat \Gr$ from $\Gr$ to the
  quotient of the free category on the directed graph $\widetilde \Gr$
  by the equivalence relation $\sim$.
  The category of partial functors $F : \Gr \to \D$
  from the partial category $\Gr$ to the partial category $\D$ is
  isomorphic to the category of 
  functors $\widehat F : \widehat \Gr \to \widehat \D$. More precisely,
  given a partial functor $F : \Gr \to \D$
  there exists a unique functor $\widehat F : \widehat \Gr \to
  \widehat \D$ such
  that $\widehat F \circ i_{\Gr} = i_{\D} \circ F$.  
\end{remark}
\begin{example}
  Given a category $\D$ there are many partial categories contained in it.
  For example if $F : \C \to \D$ is a functor then $F(\C)
  \subseteq \D$ is a partial category, but it is not in general 
  a category.
\end{example}
In order to use partial categories to codify coherence we need to
consider partial pseudo-functors. The
definition below is a variation of Grothendieck's concept of 
pseudo-functors as presented for example in \cite[Section 7.5]{Borceux}.
\begin{defn}
  A partial pseudo-functor $F$ defined on a partial
  category $\D$ consists of a function $F : \ob \D \to \ob \Cat$, a
  function $F: \mor \D \to \mor \Cat$, for composable arrows 
  $a \xrightarrow f
  b \xrightarrow g c$ of $\D$ a natural isomorphism $\gamma_{g,f} :
  F(f) \circ F(g) \to F(f\circ g)$ and for every object $a$ of $\D$ a
  natural isomorphism $\delta_a : 1_{Fa} \to F(1_a)$. These natural
  isomorphisms are required to satisfy the following coherence axioms
  for associativity and unit:
  \begin{description}
  \item[Associativity] for every triple $a \xrightarrow f b
  \xrightarrow g c \xrightarrow h d$ of composable arrows in $\D$ the
  diagram: 
  \begin{displaymath}
  \begin{CD}
    Fh \circ Fg \circ Ff @>{\id_{Fh} * \gamma_{g,f}}>> Fh \circ F(gf)
    \\
    @V{\gamma_{h,g} * \id_{Ff}}VV @V{\gamma_{h,gf}}VV \\
    F(hg) \circ Ff @>{\gamma_{hg,f}}>> F(hgf)
  \end{CD}
  \end{displaymath}
  commutes.
  \item[Unit] for every $f : a \to b$ in $\D$ the diagrams
  \begin{displaymath}
  \begin{CD}
    Ff \circ 1_{Fa} @>{\id_{Ff}* \delta_a}>> Ff \circ F1_a \\
    @V{\id_{Ff}}VV @V{\gamma_{f,1_a}}VV \\
    Ff @>{\id_{Ff}}>> F(f \circ 1_a),
  \end{CD}
  \qquad
  \begin{CD}
    1_{Fb} \circ Ff @>{\delta_b *\id_{Ff}}>> F1_b \circ Ff \\
    @V{\id_{Ff}}VV @V{\gamma_{1_b,f}}VV \\
    Ff @>{\id_{Ff}}>> F(1_b \circ f)
  \end{CD}
  \end{displaymath}
  commute.    
  \end{description}
\end{defn}
The following concept of simple morphisms was suggested by Steiner in
\cite{Steiner}.
\begin{defn}
  A
  morphism $\phi = [X \xleftarrow {\phi_1} A \xrightarrow {\phi_2} B
  \xrightarrow 
  {\phi_3} Y]$ in $U^G_+$ is called {\em simple} if  for every $y \in Y$
  the polynomial $\phi_y = \sum_{b \in {\phi_3}^{-1}(y)} (\prod_{a \in
  {\phi_2}^{-1}(b)} \phi_1(a)) \in \Z[X]$ is {\em simple}, 
  that is, 
  a sum of distinct square-free monomials.
  We denote by $sU^{G}_+$ the partial category with morphisms
  given by the class of simple morphisms in 
  $U^{G}_+$. 
\end{defn}
\begin{lem}
\label{simpleenough}
  Given a product-preserving set-valued partial functor $F$ on
  $sU^G_+$ there exists a unique extension of $F$ to a
  semi Tambara functor.
\end{lem}
\begin{proof}
  By \cite[Proposition 7.3]{Tambara} the category $U^G_+$ has a
  presentation only involving simple morphisms.
\end{proof}
\begin{defn}
\label{Tamcat}
  A {\em Tambara category} is a
  product-preserving partial
  pseudo-func\-tor on $sU^{G}_+$.
\end{defn}
We shall mostly be interested in Tambara categories arising as
homotopy categories of Quillen model categories. 
Given a Quillen model category $\D$ there 
exists 
a localization $\ho \D$
of $\D$
with respect to the class of weak equivalences. The category $\ho \D$
is the {\em homotopy category} of $\D$.
Given a functor $F \colon \D \to \mathcal E$ from a Quillen model
category to a category $\mathcal E$, a total
left derived functor of $F$ consists of a functor $LF \colon \ho
\D \to \mathcal E$ and a natural transformation $t \colon LF \circ
\gamma \to F$ with the following universal property: for every
functor $G \colon \ho 
\D \to \mathcal E$ and natural transformation $s \colon G \circ
\gamma \to F$ there exists a unique natural transformation $s'
\colon G \to LF$ such that $s = t \circ (s' *
\gamma)$.

The
two following results will be used to construct such homotopy Tambara
categories. 
\begin{prop}
\label{derivedtambara}
  Let $\C$ be a Tambara category with a Quillen model 
  structure on $\C(X)$ for every object $X$ in $U^G_+$.
  Suppose 
  \begin{enumerate}
    \item
    for every $\phi \in sU^G_+$ the functor $\C(\phi)$ has a total
    left derived functor and
    \item
    for composable morphisms $\phi$ and $\psi$
  in $sU^G_+$ the composite of the total left derived functors of
  $\C(\phi)$ and 
  $\C(\psi)$ 
  is a total left derived functor of $\C(\psi) \circ \C(\phi)$.
  \end{enumerate}
  Then there
  is a 
  Tambara category $\ho \C$ with $\ho \C(\phi)$ defined to be a
  (chosen) total left derived functor of the functor $\C(\phi)$.  
\end{prop}
\begin{proof}
  Condition (ii) implies that for composable arrows $X \xrightarrow \phi Y
  \xrightarrow \psi Z$ of $sU^G_+$ there is an isomorphism 
  $L\C(\psi) \circ L\C(\phi) \cong L(\C(\psi) \circ \C(\phi)) \cong 
  L(\C(\psi \circ \phi))$.
  Since any two total left derived functors of $\C(\psi \circ \phi)$ are
  isomorphic by a unique isomorphism, the axioms for a Tambara
  category are readily verified.  
\end{proof}
\begin{cor}
\label{browntambara}
  Let $\C$ be a Tambara category with a Quillen model categories 
  structure on $\C(X)$ for every object $X$ in $U^G_+$.
  Suppose that for every $\phi \in sU^G_+(X,Y)$ the functor
  $\C(\phi)$ preserves (acyclic) cofibrations between cofibrant
  objects. 
  Then there is a
  Tambara category $\ho \C$ with $\ho \C(\phi)$ defined to be a
  (chosen) total left derived functor of the functor $\C(\phi)$.  
\end{cor}
\begin{proof}
  K. Brown's lemma (see e.g. \cite[Lemma 9.9]{Dwyer-Spalinski}) implies that
  $\C(\phi)$ preserves weak 
  equivalences between cofibrant objects. Therefore (see
  e.g. \cite[Proposition 9.3.]{Dwyer-Spalinski}) $\C(\phi)$ has a
  total left derived functor $L\C(\phi)$ and for every cofibrant
  object $c$ of $\C(X)$ the map $L\C(\phi)(c) \to \C(\phi)(c)$ is an
  isomorphism in $\ho \C(Y)$. Since $\C(\phi)$ preserves cofibrant
  objects it follows that, given $\psi \in sU^G_+$ such that $\phi$
  and $\psi$ are composable, the composite of the total left derived
  functors of $\C(\phi)$ and $\C(\psi)$ is a total left derived
  functor of $\C(\psi) \circ \C(\phi)$.
\end{proof}
Given a $G$-set $X$ we consider $X$ as the object set of the
{\em translation category}
$X$ with $X(x,y) = \{g \in G  \colon  gx =
y \}$. Composition in $X$ is given by multiplication in $G$. 

The following construction is of fundamental importance for us. It is related
to, and inspired by, Even's construction of 
multiplicative induction in \cite{Evens}. Greenlees and May made a
similar construction in \cite{Greenlees-May}. Our construction mainly
differs from the previous ones by the fact that we do not work with 
wreath products. 
Given a set $Y$, the {\em free $\{+,\cdot\}$-algebra} on $Y$ is the set
$\coprod_{k \ge 1} \underline A(Y)_k$ where $\underline
A(Y)_1 = Y$ and 
\begin{eqnarray*}
\underline A(Y)_{k+1} &=& 
\{ w_1 + w_2 \colon w_i \in \underline A(Y)_{k_i} \text{ and } k_1 + k_2 = k\} 
\amalg
\\
&&
\{ w_1 \cdot w_2 \colon w_i \in \underline A(Y)_{k_i} \text{ and } k_1 + k_2 = k\}.
\end{eqnarray*}
\begin{construction}
\label{theTamcat}
  Given a bimonoidal
  category $(\C_0, \Box, \diamond, n_\Box, 
  n_\diamond)$ in the sense of Laplaza \cite{Laplaza} we  construct a
  Tambara category $\C = 
  \C(\C_0)$. Here  $\Box$ is the 
  additive operation and $\diamond$ is the multiplicative operation
  so one of the isomorphisms for distributivity takes the form $c_1
  \diamond (c_2 \Box c_3) \cong (c_1 \diamond c_2) \Box (c_1 \diamond
  c_3)$. 
  Recall that $\Box$ and $\diamond$ are functors from $\C_0 \times
  \C_0$ to $\C_0$ and that $n_\Box$ and $n_\diamond$ are functors from
  the trivial category $*$ to $C_0$.

  Given a finite $G$-set $X$ we let $\C(X)$ denote the category of functors
  from the 
  translation category of $X$ to $\C_0$.
  The function $X \mapsto \C(X)$ from the set of objects of $U^G_+$ to
  the class of categories clearly preserves products. Below we
  define $\C$ on morphisms.

  We
  follow Laplaza and
  let $\underline A \{X\}$ denote
  the free $\{+,\cdot\}$-algebra on
  $X \amalg \{ n_+,n_\cdot \}$.
  Given $g \in G$ we denote by $g_*$ the endomorphism of  $\underline
  A\{X\}$ induced by the action $g \cdot \colon X \to X$.

  The category $\Fun(\C(X),\C_0)$ of functors
  from $\C(X)$ to $\C_0$ has structure of a bimonoidal category with
  operations defined pointwise. There is function $\ev$ from
  $\underline A \{X\}$ to the set of objects of $\Fun(\C(X),\C_0)$ defined
  by letting $\ev_{n_+}
  = n_\Box$, $\ev_{n_\cdot} = n_\diamond$ and $\ev_x(\alpha) =
  \alpha(x)$ for $x 
  \in X$ and $\alpha \in \C(X)$ and by requiring that
  $\ev_{w_1 \cdot w_2} =
  \ev_{w_1} \diamond \ev_{w_2}$ and 
  $\ev_{w_1 + w_2} =
  \ev_{w_1} \Box \ev_{w_2}$ for $w_1,w_2 \in \underline A \{X\}$. 
  We also consider the natural
  epimorphism $\supp : \underline A \{ X\} \to \Z[X]$ defined by letting
  $\supp(n_+) = 0$, $\supp(n_\cdot) = 1$ and
  $\supp(x) = x$ for $x \in X$ and by requiring that $\supp(w_1+w_2) =
  \supp(w_1) + \supp(w_2)$ and $\supp(w_1 \cdot w_2) =
  \supp(w_1) \cdot \supp(w_2)$ 
  for $w_1,w_2 \in \underline A \{X\}$. 

  The coherence theorem of Laplaza \cite[Section
  7]{Laplaza} implies that if $\supp(w) = \supp(w') \in \Z[X]$ is
  simple, then there exists a preferred natural
  transformation $\beta_{w,w'} : \ev_w \to \ev_{w'}$,
  and if further $\supp(w'') = \supp(w')$ then $\beta_{w',w''} 
  \circ \beta_{w,w'} = \beta_{w,w''}$. 

  We choose
  for every morphism  
  $\phi = [X \xleftarrow {\phi_1} A \xrightarrow {\phi_2} B \xrightarrow
  {\phi_3} Y]$ in $sU^G_+(X,Y)$ 
  and  every $y \in Y$ an element $w_{\phi,y} \in \underline A\{
  X\}$ with
  $\supp(w_{\phi,y}) = \phi_y$.
  If $\phi$ is an identity morphism we insist on choosing $w_{\phi,y} = y$
  for $y \in Y$. 
  Note that $\supp(g_*(w_{\phi,y})) = \supp(w_{\phi,gy})$.
  Given $\alpha \in \C(X)$, the maps 
  $\alpha(g): \alpha(x) \to \alpha(g x)$ for $x \in X$
  induce a map $g_*(\alpha) : \ev_{w_{\phi,y}}(\alpha) \to
  \ev_{g_*(w_{\phi,y})}(\alpha)$.

  We define $\C(\phi) \colon \C(X) \to \C(Y)$ 
  on objects by constructing 
  $\C(\phi)(\alpha) : Y
  \to \C_0$ for $\alpha \in \C(X)$ as follows:
  let $\C(\phi)(\alpha)(y) =
  \ev_{w_{\phi,y}}(\alpha)$ for $y \in Y$
  and define
  $\C(\phi)(\alpha)(g) : \C(\phi)(\alpha)(y) \to
  \C(\phi)(\alpha)(gy)$ as the composition 
  $$\ev_{w_{\phi,y}}(\alpha)
  \xrightarrow {g_*(\alpha)}
  \ev_{g_*(w_{\phi,y})} (\alpha)
  \xrightarrow{\beta_{g_*(w_{\phi,y}),w_{\phi,gy}}} 
  \ev_{w_{\phi,gy}}(\alpha)$$ 
  for $y \in Y$ and $g \in G$. 
  We define $\C(\phi)$ on a morphism $t \colon \alpha \to \alpha'$ in
  $\C(X)$ as follows: the maps 
  $t_x \colon \alpha(x) \to \alpha'(x)$ for $x \in X$ induce maps
  $\ev_{w_{\phi,y}} (t) \colon \ev_{w_{\phi,y}}(\alpha) \to
  \ev_{w_{\phi,y}}(\alpha')$ for $y \in Y$, and these maps define the
  morphism $\C(\phi)(t) \colon \C(\phi)(\alpha) \to \C(\phi)(\alpha')$. 

  Given morphisms $\phi \in U^{G}_+(X,Y)$ and $\psi \in
  U^{G}_+(Y,Z)$ we have that $(\psi \circ \phi)_z = \psi_z(\phi_y \, | \, y
  \in Y)$, that is, $(\psi \phi)_z$ is obtained by substituting the
  variables $y \in Y$ by the polynomials $\phi_y$ in the polynomial
  $\psi_z$. Thus $\supp(w_{\psi \circ \phi,z}) =
  \supp(w_{\psi, z}(w_{\phi,y} \, 
  | \, y \in Y))$ provided that  $\phi,\psi$ and $\psi \circ \phi$ are
  simple. We define $\gamma_{\psi,\phi} : \C(\psi) \circ \C(\phi) 
  \to \C(\psi \circ \phi)$ by $(\gamma_{\psi,\phi})_z = \beta_{w_{\psi,
  z}(w_{\phi,y} \, | \, y \in Y),w_{\psi \circ \phi,z}}$ for $z \in Z$.
  Since $\C(1_Y) = 1_{\C(Y)}$ we define $\delta_Y \colon \id_{\C(Y)}
  \to \C(\id_Y)$  to be the identity.  
  The coherence axioms for $\C$ follow directly from 
  Laplaza's coherence result.
\end{construction}
It will be convenient to use the symbol
$\C$ both for the bimonoidal category $(\C_0, \Box, \diamond, n_\Box, 
  n_\diamond)$ and the Tambara category $\C$ of \ref{theTamcat}.
For example shall consider the Tambara
category $\Topp$ of pointed topological spaces with $\Topp(X) =
\Fun(X,\Topp)$. 
\begin{defn}
\label{partfib}
  A partial fibration $\GR : \Gr \to \D$ consists of a partial functor
  $\GR$ subject to the following axioms:
  \begin{description}
  \item[Composability] If the domain of $f \in \Gr$ agrees with
  the codomain of $g \in \Gr$ then $f$ and $g$ of $\Gr$ are
  composable if and only if $\GR(f)$ and $\GR(g)$ are composable
  in $\D$. 
  \item[Fibration criterion] For every arrow $\alpha : b \to c$ in
        $\D$ and object $z$ in $\Gr$ with $\GR(z) = c$ there exists
        an arrow $f : y \to z$ in $\Gr$ such that $\GR(f) = \alpha$
        and with the property that  given an arrow $g : x \to z$ in $\Gr$
        with $\Gr(g) 
              = \alpha \circ \beta$ for composable arrows $\alpha,
              \beta$ in $\D$, there exists a unique arrow $h : x \to
              y$ in $\Gr$ such that $\GR(h) = \beta$ and $g = f \circ h$.
  \end{description}
\end{defn}
The Composability axiom implies that given a partial fibration $\GR :
\Gr \to \D$ and an object $d$ of $\D$, the fiber $\GR^{-1}(d)$, that
is, the partial category 
with object set $\{x \in \Gr \colon \GR(x) = d\}$ and morphism set
$\{\alpha \in \Gr \colon \GR(\alpha) = \id_d\}$,
are categories, not only partial categories. The fibration criterion
is the usual criterion for categorical fibrations in the sense of
Gothendieck, see e.g. \cite[Chapter 8]{Borceux2}.
The following construction is well-known, see for example
loc. cit. 
\begin{construction}
\label{constructfibration}
  Given a partial pseudo-functor $F$ on
  $\D$ we 
  construct a partial fibration $\GR(F) : \Gr(F) \to \D^{\op}$ whose
  fiber at $d \in \D$ is precisely the category $F(d)^{\op}$. 
  \begin{itemize}
  \item An object of
        $\Gr(F)$ is a pair $(d,x)$ where $d \in \D$ and $x \in F(d)$
        are respectively objects of $\D$ and $F(d)$;
  \item an arrow $(d,x) \to (d',x')$ in $\Gr(F)$ is a pair
  $(\alpha,f)$ where $\alpha : d' \to d$ and $f : F(\alpha)(x') \to x$
  are respectively arrows of $\D$ and $F(d)$;
  \item $\GR(F) : \Gr(F) \to \D^{\op}$ is the first component projection,
  thus $\GR(F)(d,x) = d$ and $\GR(F)(\alpha,f) = \alpha$.
  \end{itemize}
  We must explain how to provide $\Gr(F)$ with a partial
  category-structure. 
 Consider arrows $(\alpha,f) : (d'',x'') \to (d',x')$ and $(\beta,g) :
(d',x') \to (d,x)$ in $\Gr(F)$. Provided $\alpha$ and $\beta$
are composable this yields the composite
\begin{displaymath}
  F(\alpha \beta)(x'') {\cong} F(\beta)F(\alpha)(x'')  
  \xrightarrow 
  {F(\beta)(f)} F(\beta)(x') \xrightarrow g 
  x 
\end{displaymath}
in $F(d)$, where the first isomorphism is the associativity
isomorphism for the partial pseudo-functor $F$. Writing $g\star f$ for this
composite, 
the composition law is defined by the relation $(\alpha,f) \circ
(\beta,g) = (\beta \circ \alpha, g \star f)$. The associativity of this
composition follows immediately from the associativity axiom for 
pseudo-functors. On the other hand the unit axiom for 
pseudo-functors 
implies that the pair $(1_d,\delta_{d,x}^{-1})$ is an identity
morphism on the object $(d,x)$ of $\Gr(F)$. Here $\delta_d \colon 1_{Fd}
\to F(1_d)$ is the unit isomorphism for $F$. 
This proves that
$\Gr(F)$ is a partial category. The functoriality of
$\GR(F)$ is obvious.
\end{construction} 
We leave the proof of the following lemma to the reader.
\begin{lem}
\label{Grprod}
  Given a product-preserving partial pseudo-functor $F : \D \to
  \Cat$, objects $d,d'$ of $\D$ and an object $z$ of $F(d \times d')$,
  the object $(d \times d', z \in F(d \times d'))$ is the product of
  $(d,F(\pr_1)z)$ and $(d',F(\pr_2)z)$ in $\Gr(F)^{\op}$.
\end{lem}

\section{Induction of Cofibrations}
\label{indofcof}
In this section we consider 
a
complete and cocomplete
symmetric monoidal category  
$(\C_0,\diamond,u_\diamond)$, with the property that the functor $c
\diamond -$ 
commutes with colimits for every object $c$ of $C_0$.
Denoting the coproduct of $\C_0$ by $\amalg$ and the initial
object by $\emptyset$, we have a
bimonoidal category
$(\C_0,\amalg,\diamond,\emptyset,u_\diamond)$  and we can consider the
Tambara category  
$\C = \C(\C_0)$ of Construction \ref{theTamcat}.
Since $\diamond$ preserves colimits we have
for every finite $G$-set $X$ and every $c \in \ob \C(X)$ 
that the functor $c \diamond -$ preserves push-outs and finite coproducts
and that $c \diamond \emptyset \cong \emptyset$. 

Given a map $f: X \to Y$ of finite $G$-sets we obtain
functors $f_\diamond := C(N_f)$ and $f_\amalg := C(T_f)$ from $C(X)$
to $C(Y)$ 
and $f^* := C(R_f)$ from $C(Y)$ to $C(X)$. 

Further we fix for every finite $G$-set $X$ subsets $J(X)$ and $I(X)$ 
of the class of morphisms in $\C(X)$ and given objects $a,b$ of
$\C(X)$ we write $J(X)(a,b)$ and $I(X)(a,b)$ for $J(X) \cap
\C(X)(a,b)$ and $I(X) \cap \C(X)(a,b)$ respectively. These subsets are
required to 
satisfy firstly that given a sum diagram $X_1 \xrightarrow{i_1} X
\xleftarrow{i_2} X_2$ of finite $G$-sets, the isomorphism $(i_1^*,i_2^*) : \C({X})
\xrightarrow \cong \C({X_1}) \times  
\C({X_2})$ induces bijections 
\begin{eqnarray*}
  J(X)(a,b)  &\cong& J({X_1})(i_1^*a,i_1^*b) 
  \times J({X_2})(i_2^*a,i_2^*b), \\
  I({X})(a,b)&\cong& I({X_1})(i_1^*a,i_1^*b) 
  \times I({X_2})(i_2^*a,i_2^*b).
\end{eqnarray*}
Secondly we require that for every finite $G$-set $X$ the data
$(\C(X),I(X),J(X))$ 
specifies a cofibrantly generated symmetric monoidal Quillen model
category with 
generating cofibrations $I(X)$ and generating acyclic cofibrations
$J(X)$, and that for every $G$-map $f : X \to Y$, the pair 
$f_\amalg : \C(X) \rightleftarrows \C(Y) : f^*$ is a Quillen adjoint pair. 

By abuse of notation we shall let $\Z$ denote the category associated
to the underlying partially ordered set of the integers. The functor
$\Fun(\Z,-)$ 
taking a category $\D$ to the category 
$\Fun(\Z,\D)$ of filtered objects in $\D$ preserves products so
we obtain a Tambara
category $\Z\C := \Fun(\Z,\C)$. In Section \ref{secfilob} we have
collected some basic facts about filtered object. 
\begin{defn}
  Given a $G$-set $X$ and a morphism $c: C_{-1} \to C_{0}$ in $\C(X)$,
  we denote by $\overline C$ the object in $(\Z\C)(X)$ with $\overline
  C(i) = C_{-1}$ for $i 
  \le -1$ and with $\overline C (i) = C_0$ for $i \ge 0$.
  For $i=-1$ the morphism $\overline C(i,i+1) \colon \overline C(i) \to
  \overline C(i+1)$ 
  is the morphism $c$ and for $i \ne -1$ it is an identity morphism.
\end{defn}
\begin{thm}
\label{maininduction}
  Suppose that for every $G$-map $g \colon W \to Z$ the map
  $(g_{\diamond}\overline C)(-1) \to (g_{\diamond}
  \overline C)(0)$ is a
  cofibration for every
  $c$ in $\in I(W)$ and the map 
  $(g_{\diamond}\overline C)(-1) \to (g_{\diamond}
  \overline C)(0)$ is an
  acyclic cofibration for every
  $c$ in $\in I(W)$.
  Then for
  every map
  $f : X \to Y$ of finite $G$-sets the functor $f_{\diamond} : \C(X) \to \C(Y)$
  preserves both cofibrations and acyclic cofibrations between
  cofibrant objects. 
\end{thm}
\begin{cor}
\label{htytambara0}
  Suppose for every map $g
  : W \to Z$ of finite $G$-sets that 
  \begin{enumerate}
    \item
    for every $d : D_{-1} \to D_0$ in
  $I(Z)$ the map $g^*(D_{-1}) \to g^*(D_0)$ is
  an  
  cofibration in $\C(W)$, 
    \item
    for every $d : D_{-1} \to D_0$ in
  $J(Z)$ the map $g^*(D_{-1}) \to g^*(D_0)$ is
  an acyclic 
  cofibration in $\C(W)$, 
    \item
    for every $c :
  C_{-1} \to C_{0}$ in $I(W)$ 
  the map $(g_{\diamond}\overline C)(-1) \to (g_{\diamond}
  \overline C)(0)$ is a cofibration in $\C(Z)$ and
    \item
    for every $c :
  C_{-1} \to C_{0}$ in $J(W)$ 
  the map $(g_{\diamond}\overline C)(-1) \to (g_{\diamond}
  \overline C)(0)$ is an acyclic cofibration in $\C(Z)$.
  \end{enumerate}
  Then for every for every morphism $\phi \in sU^G_+$ the functor
  $\C(\phi)$ preserves both cofibrations and acyclic cofibrations
  between cofibrant 
  objects. 
\end{cor}
\begin{proof}
  For a map $f : X \to Y$ of finite $G$-sets the functor $f^* \colon
  \C(Y) \to \C(X)$ has a
  right adjoint $f_\Pi \colon \C(X) \to \C(Y)$ with $f_\Pi(C)(y) 
  \cong \prod_{x \in f^{-1}(y)} C(x)$. In particular $f^*$ preserves
  relative cell complexes. Since every (acyclic) cofibration is a
  retract of a relative cell complex the functor $f^*$ preserves
  (acyclic) 
  cofibrations. Thus  
  for every map $f : X \to Y$ of finite $G$-sets the functors
  $f_\amalg$ and $f^*$ preserve (acyclic) cofibrations between
  cofibrant objects, and so does $f_\diamond$ by \ref{maininduction}.
  By \cite[Proposition 7.3]{Tambara} for every $\phi \in sU_+^G$ the
  functor $\C(\phi)$ is isomorphic to a
  composite of functors of the above 
  sort. 
\end{proof}
\begin{cor}
\label{htytambara}
  For every symmetric monoidal category $\C_0$ with symmetric monoidal
  Quillen model structures
  on the categories $\C(X) = \Fun(X,\C_0)$ for $X$ in $U^G_+$
  satisfying the assumptions of 
  Corollary \ref{htytambara0} there exists a
  Tambara category $\ho \C$ with $\ho \C(\phi)$ defined to be a
  (chosen) total left derived functor of the functor $\C(\phi)$.  
\end{cor}
\begin{proof}
  This follows by combining Corollary \ref{htytambara0} and Corollary
  \ref{browntambara}. 
\end{proof}
Since we are only dealing with finite $G$-sets, Theorem
\ref{maininduction} is a consequence of the following theorem:
\begin{thm}
\label{maininduction2}
  Let $n \in \N$. 
  Suppose that for every $G$-map $g \colon W \to Z$ with fibers of
  coarinality les than or equal to $n$ the map
  $(g_{\diamond}\overline C)(-1) \to (g_{\diamond}
  \overline C)(0)$ is a
  cofibration for every
  $c$ in $\in I(W)$ and the map 
  $(g_{\diamond}\overline C)(-1) \to (g_{\diamond}
  \overline C)(0)$ is an
  acyclic cofibration for every
  $c$ in $\in I(W)$.
  Then for every
  map $f : X \to Y$ of finite $G$-sets with fibers of cardinality less
  than or equal to 
  $n$, the functor $f_{\diamond} : \C(X) \to \C(Y)$
  preserves both cofibrations and acyclic cofibrations between
  cofibrant objects. 
\end{thm}
\begin{example}
  We illustrate Theorem
  \ref{maininduction2} by considering a simple example.
  Let $G = \Z/2\Z$ and
  consider the map $f \colon G/e \to G/G$ of $G$-sets. The
  functor $\C_0 =  \Fun(*,\C_0) \to \Fun(G/e,\C_0) = \C(G/e)$ induced
  by the functor from the translation category of $G/e$ to the final
  category $*$ is an equivalence of categories. 
  Suppose that
  \begin{displaymath}
  \begin{CD}
    C_{-1} @>c>> C_0 \\
    @VVV @VVV \\
    A_{-1} @>a>> A_0
  \end{CD}
  \end{displaymath}
  is a push-out diagram in $\C_0$ and that the map
  $(f_\diamond \overline C)(-1) \to
    (f_\diamond \overline C)(0)$
  is a cofibration in $\C(G/G)$. Then also the lower horizontal map in
  the push-out diagram 
  \begin{displaymath}
  \begin{CD}
    (f_\diamond \overline
    C)(-1) 
    @= 
    C_{-1} \diamond C_0 \coprod_{C_{-1} \diamond C_{-1}} C_0 \diamond C_{-1} 
    @>>> (f_\diamond \overline C)(0) \\
    @VVV @. @VVV \\
    (f_\diamond \overline
    A)(-1) 
    @=
    A_{-1} \diamond A_0 \coprod_{A_{-1} \diamond A_{-1}} A_0 \diamond A_{-1} 
    @>>> (f_\diamond \overline A)(0)
  \end{CD}
  \end{displaymath}
  is a cofibration in $\C(G/G)$. On the other hand, suppose that the map
  $C_{-1} \diamond C_{-1} \to C_{0} \diamond C_{-1}$ is a cofibration
  in $\C_0$. Then $f_{\amalg}(C_{-1} \diamond C_{-1}) \to
  f_{\amalg}(C_{0} \diamond C_{-1})$ is a cofibration in $\C(G/G)$, and 
  therefore also the lower horizontal map in the 
  push-out diagram 
  \begin{displaymath}
    \begin{CD}
      f_{\amalg}(C_{-1} \diamond C_{-1}) @>>> 
      f_{\amalg}(C_{0} \diamond C_{-1}) \\
      @VVV @VVV \\
     A_{-1} \diamond A_{-1} @>>>  
     A_{-1} \diamond A_0 \coprod_{A_{-1} \diamond A_{-1}} A_0 \diamond A_{-1} 
    \end{CD}
  \end{displaymath}
  is a cofibration in $\C(G/G)$. Thus the composition
  \begin{displaymath}
    f_{\diamond}(A_{-1}) = A_{-1} \diamond A_{-1} \to     (f_\diamond
    \overline A)(-1) \to     (f_\diamond \overline
    A)(0) =     f_\diamond(A_0) 
  \end{displaymath}
  is a cofibration.
\end{example}
The rest of this section is devoted to a proof by induction
on $n$ of Theorem \ref{maininduction2}. 
\begin{lem}
\label{useretract}
  Let $f : X \to Y$ be a map of finite $G$-sets.
  Suppose that the map 
  $(f_\diamond \overline C)(-1) \to (f_\diamond \overline
  C)(0)$ is a cofibration in $\C(Y)$
  for every $c : C_{-1} \to C_0$ in $I(X)$ and that
  it is an acyclic cofibration 
  for every $c : C_{-1} \to C_0$ in  $J(X)$.
  Then the map 
  $(f_\diamond \overline
  A)(-1) \to (f_\diamond \overline A)(0)$ is a cofibration
  in $\C(Y)$ 
  for every 
  cofibration $a : A_{-1} \to A_0$ in $\C(X)$ and it is an acylcic
  cofibration 
  for every acyclic cofibration $a : A_{-1} \to A_0$ in $\C(X)$.
\end{lem}
\begin{proof}
  Assume first that there exists a push-out diagram of the form
  \begin{displaymath}
  \begin{CD}
    C_{-1} @>c>> C_0 \\
    @VVV @VVV \\
    A_{-1} @>a>> A_0
  \end{CD}
  \end{displaymath}
  with $c$ in $I(X)$ (respectively $J(X)$). Then by Lemma \ref{cralem2a} the
  diagram
  \begin{displaymath}
  \begin{CD}
    (f_\diamond \overline
    C)(-1) @>>> (f_\diamond \overline C)(0) \\
    @VVV @VVV \\
    (f_\diamond \overline
    A)(-1) @>>> (f_\diamond \overline A)(0)
  \end{CD}
  \end{displaymath}
  is a
  push-out diagram.
  Since by assumption the upper horizontal map is a cofibration, so is
  the lower horizontal map. It follows that the map $    (f_\diamond \overline
    A)(-1) \to (f_\diamond \overline A)(0)$ is an (acyclic)
  cofibration for every relative cell complex $a : A_{-1} \to
  A_0$. The lemma now follows from the fact that every (acyclic)
  cofibration is a retract of a relative cell complex.
\end{proof}
\begin{defn}
\label{hatdefn}
  Given a morphism $c: C_{-1} \to C_{0}$ in $\C(X)$, we define 
  an object $\widehat C$ in $(\Z\C)(X)$ as follows: for $i < -1$ we
  let $\widehat C(i) = \emptyset$ be the initial object of $\C(X)$, for
  $i = -1$ we let $\widehat C(i) = C_{-1}$, and for $i \ge 0$ we let
  $\widehat C (i) (x) = C_0$. The map $\widehat C(i) \to
  \widehat C(i+1)$ is the identity except for $i = -2,-1$. For $i =
  -2$ it is
  the unique map
  $\emptyset = \widehat C(-2)(x) \to \widehat C(-1)(x)$, and for $i =
  -1$ it is given by the map $c$. 
\end{defn}
Note that for every map $f:X \to Y$ of finite $G$-sets there is
an isomorphism $(f_\diamond \overline C)(-1) \cong (f_\diamond
\widehat C)(-1)$. 
\begin{lem}
\label{usepop}
  Suppose that for every map $g : W \to Z$ of finite $G$-sets
  with fibers of
  cardinality less than or equal to $n-1$ we have: 
  \begin{enumerate}
    \item
    for 
  every cofibrant object 
  $B$ of $\C(W)$ the 
  object $g_\diamond B$ is cofibrant in $\C(Z)$, 
    \item
    for
  every $c : C_{-1} \to C_0$ in $I(W)$ the map
  $(g_\diamond \overline C)(-1) \to (g_\diamond \overline C)(0)$ is a
  cofibration and 
    \item
    for
  every $c : C_{-1} \to C_0$ in $J(W)$ the map
  $(g_\diamond \overline C)(-1) \to (g_\diamond \overline C)(0)$ is an
  acyclic cofibration. 
  \end{enumerate}
  Then for every map $f : X
  \to Y$ of finite $G$-sets with fibers of cardinality $n$
  and every $k \in \Z$ 
  with $-n < k < 0$, the map $(f_\diamond \widehat A)(k-1) \to
  (f_\diamond \widehat 
  A)(k)$ is a cofibration for
  every 
  cofibration $a : A_{-1} \to A_0$ of cofibrant objects in $\C(X)$ and
  it is an acyclic cofibration for 
  every 
  acyclic cofibration $a : A_{-1} \to A_0$ of cofibrant objects in $\C(X)$.
\end{lem}
\begin{proof}[Proof of theorem \ref{maininduction2}]
  We prove the theorem by induction on $n$. The theorem holds in the
  case $n =1$. Assume that the theorem holds for $n-1$. Using the
  push-out product axiom in $\C(Y)$ we see that it
  suffices to show that it holds for every
  map $f : X \to Y$ of finite $G$-sets with fibers of cardinality
  exactly $n$
  and every (acyclic) cofibration $a : A_{-1} \to A_0$ of cofibrant
  objects in $\C(X)$. In that case
  $f_\diamond A_{-1} \cong (f_\diamond \widehat A)(-n)$ and 
  $f_\diamond A_{0} \cong (f_\diamond \widehat A)(0)$. It follows that
  it suffices to show that for $-n+1 \le k \le 0$ the map 
  $(f_\diamond \widehat A)(k-1) \to (f_\diamond \widehat A)(k)$ is
  an (acyclic) cofibration in $\C(Y)$. The case $k = 0$ is treated in
  Lemma \ref{useretract}, and using our inductive assumption, the case
  $-n < k<0$
  follows from Lemma \ref{usepop}. 
\end{proof}
We now introduce some notation needed for the proof of Lemma \ref{usepop}.
Let $f : X \to Y$ be a map of finite $G$-sets and let $p: \{-1,0\} \times X
\to X$ denote the  
projection. We consider an exponential diagram of the form:
\begin{displaymath}
\begin{CD}
  X @<p<< \{-1,0\} \times X @<e<< X'\\
  @VVfV @. @VVf'V \\
  Y @= Y @<q<< Y'.
\end{CD}
\end{displaymath}
Given $y' \in Y'$ we let $|y'| := \sum_{x'\in
(f')^{-1}(y')} \pr (e(y')) \in \Z$,
where $\pr : \{-1,0\} \times X \to \{-1,0\}$ is the 
projection, and given $k \in \Z$ we let 
$i_k : Y'_{k} \to Y'$ 
denote the inclusion of the subset $Y'_k$ consisting of the
elements $y' \in Y'$ with $|y'| = k$. 
We let $X'_{-1}$ and $X'_0$ be defined as pull-backs in the diagram
\begin{displaymath}
\begin{CD}
  X'_{-1} @>>> X' @<<< X'_0 \\
  @VVV @VVV @VVV \\
  \{-1\} \times X @>>> \{-1,0\} \times X @<<< \{ 0\} \times X.
\end{CD}
\end{displaymath} 
We can now construct a commutative diagram of the form: 
\begin{displaymath}
\begin{CD}
  X' @= X'_{-1} \amalg X'_0
  @<{j_{-1,k}\amalg j_{0,k}}<< X'_{-1,k} \amalg X'_{0,k} \\
  @VVf'V @VV{f'_{-1}\amalg f'_0}V @VV{f'_{-1,k}\amalg
  f'_{0,k}}V \\
  Y'@<{\nabla}<< Y' \amalg Y' @<{i_k \amalg i_k}<<
  Y'_k \amalg Y'_k,
\end{CD}
\end{displaymath}
where $\nabla$ is the fold map and where the right hand square is
a pull-back.

Given a map
$a: A_{-1} \to A_0$ in $\C(X)$ we use the isomorphism
$\C(X)
\times \C(X) \cong \C(\{-1\} \times X) \times \C(\{0 \} \times X)
\cong \C(\{-1,0\}\times X)$
to consider $(\emptyset,A_{-1})$ and $(A_{-1},A_{0})$ as objects of
$\C(\{-1,0\} \times X)$. We 
let $A(k)_{-1} = (j_{-1,k} \amalg j_{0,k})^* e^*
(\emptyset,A_{-1})$
and 
$A(k)_{0} = (j_{-1,k} \amalg j_{0,k})^* e^*
(A_{-1},A_{0}) \in
\C({X'}_{-1,k} \amalg {X'}_{0,k})$. The map $a$ induces a map $a(k) :
A(k)_{-1} \to A(k)_0$, and we consider the object $\widehat {A(k)} \in
\Z\C({X'}_{-1,k} \amalg {X'}_{0,k})$ of Definition \ref{hatdefn}. Note that
considering $A_{-1}$ 
as a functor from $\Z$ to $\C(X)$ with $A_{-1}(k) = \emptyset$ for $k
< 0$ and $A_{-1}(k) = A_{-1}$ for $k \ge 0$ we have that $\widehat
{A(k)} = (j_{-1,k}^* e^* p^* A_{-1}, j_{0,k}^* e^* p^* \widehat A)$.
We define an object $T_{A,k}$ of $\Z\C(Y)$ by the formula 
$$T_{A,k} := q_\amalg  \circ (\nabla \circ (i_k \amalg i_k) \circ (f'_{-1,k}\amalg
f'_{0,k}))_\diamond (\widehat{A(k)}).$$
\begin{lem}
\label{takpo}
  For every $k \in \Z$ there is a push-out diagram in $\C(Y)$ of the
  form 
  \begin{displaymath}
  \begin{CD}
    T_{A,k}(-1) @>>> T_{A,k}(0) \\
    @VVV @VVV \\
    (f_\diamond \widehat A)(k-1) @>>> (f_\diamond \widehat A)(k).
  \end{CD}
  \end{displaymath}
\end{lem}
\begin{proof}
  Fix $y \in Y$ and let $U = 
  f^{-1}(y)$.
  There is a
  bijection $\phi : \{-1,0\}^{U} \xrightarrow \cong 
  q^{-1}(y)$. Let $V_k \subseteq V_{\le k}
  \subseteq \Z^U$ denote the sets consisting of the maps $a \colon U \to \Z$
  with $\sum_{u \in U} a(u) = k$ and $\sum_{u \in U} a(u) \le k$ respectively.
  There is a unique functor $T : (\Z^U,\le) \to \C(*)$ with 
  $T(\alpha) = \underset {x \in U} \Diamond A_{\alpha(x)}(x)$. 
  For $\alpha : U \to \{-1,0\}$ we have
  \begin{displaymath}
   T(\alpha) 
   =
   \underset {x \in U} \Diamond A_{\alpha(x)}(x) 
   \cong \underset {x' \in {f'}^{-1}(\phi(\alpha))} \Diamond
   (e^*(A_{-1},A_0))(x')  
   \cong 
   (f'_{\diamond} e^* 
   (A_{-1},A_{0}))(\phi(\alpha)). 
  \end{displaymath}
  Writing out the definitions we get 
  \begin{eqnarray*}
    (T_{A,k}(0))(y) &\cong& \coprod_{y'\in q^{-1}(y)} ((\nabla(i_k
    \amalg i_k) (f'_{-1,k} \amalg f'_{0,k}))_{\diamond} \widehat{A(k)}
    )(0)(y') \\
    &\cong&
    \coprod_{y'\in q^{-1}(y)} (i_k f'_{-1,k})_\diamond
    (j^*_{-1,k}e^*p^* A_{-1})(y') 
    \diamond
    (i_k f'_{0,k})_\diamond
    (j^*_{0,k}e^*p^* \widehat A)(0)(y') \\
    &\cong&
    \coprod_{\alpha \in V_k \cap \{-1,0\}^U} 
    (\underset {x \in \alpha^{-1}(-1)} \Diamond A_{-1}(x)) 
    \diamond
    (\underset {x \in \alpha^{-1}(0)} \Diamond A_{0}(x)) \\     
    &\cong&
    \coprod_{\alpha \in V_k \cap \{-1,0\}^U} 
    (\underset {x \in U} \Diamond A_{\alpha(x)}(x)) 
    = 
    \coprod_{\alpha \in V_k \cap \{-1,0\}^U} 
    T(\alpha),
  \end{eqnarray*}  
  and similarly we get
  \begin{eqnarray*}
    (T_{A,k}(-1))(y) &\cong& 
    \coprod_{\alpha \in V_k \cap \{-1,0\}^U} 
    (\underset {x \in \alpha^{-1}(-1)} \Diamond A_{-1}(x)) 
    \diamond
    ((\underset {x \in \alpha^{-1}(0)} \Diamond \widehat A
    (x))(-1)) \\      
    &\cong&
    \coprod_{\alpha \in V_k \cap \{-1,0\}^U} 
    \colim{\beta \in V_{< \alpha}} T(\beta),
  \end{eqnarray*}  
where $V_{< \alpha} \subseteq \Z^U$ consists of those $\beta \in \Z^U
\setminus \alpha$ satisfying $\beta(u) \le \alpha(u)$ for $u \in U$.
Note that if $\alpha \notin \{-1,0\}^U$ then the natural map
$\colim{\beta \in V_{< \alpha}} T(\beta) \to T(\alpha)$ is an isomorphism.
On the other hand we have  
\begin{eqnarray*}
  (f_\diamond \widehat A)(k)(y) &\cong& (\underset {x \in U}
  \Diamond \widehat A(x))(k) 
  \cong \colim{\alpha \in V_{\le k} \cap \{-1,0\}^U} (\underset {x \in U}
  \Diamond 
  A_{\alpha(x)}(x)) \cong \colim{\alpha \in V_{\le k}} T(\alpha),
\end{eqnarray*}
and $
  (f_\diamond \widehat A)(k-1)(y) \cong \colim{\beta \in V_{\le
  {k-1}}} T(\beta).
$
The lemma now follows from Lemma \ref{pufilt}.
\end{proof}
\begin{lem}
\label{takcof}
  Under the assumptions of Lemma \ref{usepop} the map 
  $T_{A,k}(-1) \to T_{A,k}(0)$ is an (acyclic) cofibration for $-n < k
  < 0$. 
\end{lem}
\begin{proof}
  Since $-n<k<0$ the fibers of the maps $f'_{-1,k}$ and $f'_{0,k}$ are of
  cardinality less than or equal to
  $n-1$. By assumption the object $(f'_{-1,k})_\diamond j_{-1,k}^* e^*
  p^* A_{-1}$ is 
  cofibrant in $\C(Y'_k)$. The map $a(k) = (j_{-1,k} \amalg
  j_{0,k})^* e^* (\id_{A_{-1}},a)$ 
  is an (acyclic) cofibration 
  and so is by
  Lemma \ref{useretract} the map 
  $$(f'_{-1,k} \amalg f'_{0,k})_\diamond \widehat {A(k)}(-1) \to 
  (f'_{-1,k} \amalg f'_{0,k})_\diamond \widehat {A(k)}(0)$$
  in $\C(Y'_k)$.
  Under the isomorphism 
  $\C(Y'_k) \times \C(Y'_k) \cong \C(Y'_k \amalg
  Y'_k)$
  this map can be identified with the map
  \begin{eqnarray*}
  && ((f'_{-1,k})_\diamond j_{-1,k}^* e^*
  p^* A_{-1},((f'_{0,k})_\diamond j_{-1,k}^* e^*
  p^* \widehat
  A)(-1)) \to \\
  && ((f'_{-1,k})_\diamond j_{-1,k}^* e^*
  p^* A_{-1},((f'_{0,k})_\diamond j_{-1,k}^* e^*
  p^* \widehat
  A)(0)).   
  \end{eqnarray*}
  and 
  the map
  \begin{displaymath}
    (\nabla (i_k\amalg i_k) (f'_{-1,k} \amalg f'_{0,k}))_\diamond
    \widehat {A(k)} (-1)
    \to
    (\nabla (i_k\amalg i_k) (f'_{-1,k} \amalg f'_{0,k}))_\diamond
    \widehat {A(k)} (0)
  \end{displaymath}
  can be identified with the map
  \begin{eqnarray*}
  && (i_k \circ f'_{-1,k})_\diamond j_{-1,k}^* e^*
  p^* A_{-1} \diamond ((i_k \circ f'_{0,k})_\diamond j_{0,k}^* e^*
  p^* \widehat
  A )(-1) \to \\
  && (i_k \circ f'_{-1,k})_\diamond j_{-1,k}^* e^*
  p^* A_{-1} \diamond ((i_k \circ f'_{0,k})_\diamond j_{0,k}^* e^*
  p^* \widehat
  A)(0).
  \end{eqnarray*}
  This map is an (acyclic) cofibration by the push-out-product
  axiom in $\C(Y')$. The result now follows from the fact that
  $q_\amalg $ being a left 
  Quillen functor it preserves (acyclic) cofibrations.
\end{proof}
\begin{proof}[Proof of Lemma \ref{usepop}]
  By \ref{takpo} and \ref{takcof}, for $-n < k < 0$ the map
  $(f_\diamond \widehat 
  A)(k-1) \to (f_\diamond \widehat
  A)(k)$ is a push-out of an (acyclic) cofibration. The statement of
  Lemma \ref{usepop} now follows since a push-out of an (acyclic)
  cofibration is an (acyclic) cofibration.
\end{proof}
\section{Constructing Tambara Functors}
\label{firsttrfsec}
In this section we study 
(semi-) Tambara 
functors arising from Tambara categories. 
Given a Tambara category $\C$ there is an opposite Tambara functor
$\C^{\op}$ with $\C^{\op}(X) = \C(X)^{\op}$. 
In Section \ref{tamcatsec} we constructed
product-preserving 
partial functors $\GR(\C^{\op})^{\op} : \Gr(\C^{\op})^{\op} \to sU^G_+$
and  $\GR(\C)^{\op} : \Gr(\C)^{\op} \to sU^G_+$. 
We denote by 
$\hom_\C : (\Gr(\C^{\op}) \times_{{U^G_+}^{\op}}
\Gr(\C))^{\op} \to \Set$ 
the product-preserving partially
defined functor 
with value
$\hom_\C((X,x),(X,y)) = \C(X)(x,y)$ on the object $((X,x),(X,y))$ of
$(\Gr(\C^{\op}) \times_{{sU^G_+}^{\op}} 
\Gr(\C))^{\op}$. 
Let $s : sU^G_+ \to
(\Gr(\C^{\op}) \times_{{sU^G_+}^{\op}} \Gr(\C))^{\op}$ be a
product-preserving section of the  
projection $(\Gr(\C^{\op}) \times_{{sU^G_+}^{\op}} \Gr(\C))^{\op} \to
sU^G_+$.  
By Lemma \ref{simpleenough} 
there exists a unique
semi Tambara functor ${\C(s)} : U^G_+
\to \Set$ extending the
product-preserving partial functor 
$\hom_\C \circ  s : sU^G_+  \to \Set$.
\begin{defn}
A {\em commutative semi-ring} in a bimonoidal category $\C_0$
is an object $A$ of 
$(\C_0,\diamond,\Box,n_\diamond,n_\Box)$ with
commutative monoid structures for both $\diamond$ and $\Box$
making the following distributivity diagram commutative:
\begin{displaymath}
\begin{CD}
  A \diamond (A \Box A) @>>> A \diamond A @>>> A\\
  @V{\cong}VV @. @| \\
  (A \diamond A) \Box (A \diamond A) @>>> A \Box A @>>> A.
\end{CD}
\end{displaymath}  
\end{defn}

\begin{construction}
\label{commonsection}
Let $(\C_0,\diamond,\Box,n_\diamond,n_\Box)$ be
a bimonoidal category, let $\C = \C(\C_0)$ denote the Tambara
category of \ref{theTamcat}, let $A$ be a commutative semi-ring in
$\C_0$ and let $B$ be a commutative semi-ring $\C_0^{\op}$.
We shall construct a product-preserving partial 
functor 
$$(s_B,s_A) : sU^{G}_+
\to (\Gr(\C^{\op}) 
\times_{(sU^{G}_+)^{\op}} \Gr(\C))^{\op}.$$
As in the above discussion we get a 
product-preserving partial functor  
\begin{displaymath}
  sU^{G}_+ \xrightarrow {(s_B,s_A)}  (\Gr(\C^{\op})
  \times_{(sU^{G}_+)^{\op}} \Gr(\C))^{\op} 
  \xrightarrow {\hom_{\C}} \Set
\end{displaymath}
and a semi Tambara 
functor ${\C}(B,A) : U^{G}_+ \to
\Set$ extending the above partial functor.
Note that ${\C}(A,B)(X) = \C(X)(p^*_XB,p^*_XA)$.

Here $p_X : X \to *$ denotes the $G$-map from $X$ to a point
and $p_X^*$ denotes the functor $\C(R_{p_X}) : \C(*) \to \C(X)$. 
We define a section 
$$s_A : 
sU^{G}_+ \to \Gr(\C)^{\op}$$ 
of 
$\GR(\C)^{\op}$ with $s_A(X) = 
(X,p_X^*A)$.
By Lemma \ref{Grprod}
$s_A$ preserves products. 
Let $\phi : X \to Y$ be a morphism of
$sU^{G}_+$.
In order to construct $s_A(\phi) \colon (Y,p_Y^*A) \to (X,p_X^*A)$ we
need to specify a morphism $s_A(\phi)(y) \colon \C(\phi)(p_X^*A)(y) 
\to (p_Y^*)(y) = A$ for every $y \in Y$. 
In \ref{theTamcat}  we 
chose $w_{\phi,y} \in
\underline A \{ X\}$  
with $\supp(w_{\phi,y}) = \phi_y$ and we defined the element 
$\C(\phi)(p_X^*A)(y) =
\ev_{w_{\phi,y}} (p^*_X)$ as the evaluation in $\C_0$ of the word
$w_{\phi,y}$ in
the letters $A,\Box,\diamond, n_\Box$ and $n_\diamond$.
Since $A$ is a
a monoid with respect to both $\Box$ and $\diamond$ there is a
morphism $s_A(\phi)(y) \colon \ev_{w_{\phi,y}} (p^*_X) \to A$, and
this morphism is 
independent of the choice of $w_{\phi,y}$ because $A$ is a commutative
semi-ring. This ends the construction of $s_A$.

Dually, since $B$ is a commutative semi-ring in 
$\C_0^{\op}$, we obtain a product-preserving section 
$s_B : 
sU^{G}_+ \to {\Gr(\C^{\op})}^{\op}$ of 
$\GR(\C^{\op})^{\op}$.
Combining $s_A$ and $s_B$ we obtain the desired product-preserving partial 
functor $(s_B,s_A) : sU^{G}_+
\to (\Gr(\C^{\op}) 
\times_{(sU^{G}_+)^{\op}} \Gr(\C))^{\op}$. 
\end{construction}
\begin{example}
\label{invariant-ring}
The category
  $(\Ab,\oplus,\otimes,0,\Z)$ of 
abelian groups
is 
bimonoidal.
A semi-ring $A$ in $\Ab$ is the same
thing as a commutative algebra, and a semi-ring
$B$ in $\Ab^\op$ is the same thing as a commutative coalgebra. For a finite
$G$-set of the form 
$G/H$ for a subgroup $H$ of  $G$ we have
that $\hom_{\Ab}(B,A)(G/H) = \hom_{\Z}(B,A)^{H}$ is the ring of
$H$-equivariant $\Z$-linear maps from $B$ to $A$. In particular for $B
= \Z$ we recover the invariant ring Tambara functor of \cite[Example
  3.1]{Tambara} 
\end{example}
\begin{example}
\label{evens-example1}
 A semi-ring in the category $(\Ch_{\Z},\oplus,\otimes,0,\Z)$ of chain
complexes of abelian
groups is the same
thing as a differential graded commutative algebra, and a semi-ring in
$\Ch^\op_\Z$ is the same thing as
a differential graded commutative coalgebra. For a
finite $G$-set of the form 
$G/H$ for a subgroup $H$ of  $G$ we have
that $\hom_{\Ch_\Z}(B,A)(G/H) = \Ch_{\Z}(B,A)^{H}$ is the ring of
$H$-equivariant $\Z$-linear chain maps from $B$ to $A$. 
It is desirable to pass to homology of these chain complexes. 
Evens does this in the case 
$B = \Z[T]$ where $T$ in degree $2$ has trivial action of $G$, and
where $A$ is concentrated in degree zero \cite{Evens},
\cite[Example 3.4]{Tambara}.  
We shall return to this situation in example \ref{evens-example}.
\end{example}
Given a Tambara
category $\C$ satisfying the assumptions of Proposition
\ref{derivedtambara} 
there exists a Tambara category $\ho \C$ with $\ho \C(\phi)$ given by a total
left derived functor of $\C(\phi)$ for $\phi \in
sU^G_+$.
\begin{prop}
\label{monoidtam}
  Under the assumptions of Proposition
  \ref{derivedtambara} there is a partial functor $\gamma : \Gr(\C)
  \to \Gr(\ho \C)$ 
  with $\GR(\ho \C) \circ \gamma = \GR(\C)$.
\end{prop}
\begin{proof}
  On objects we define $\gamma$ to be the identity, and we define
  $\gamma$ on morphisms by letting $\gamma$ take a
  morphism $(X,c) \xleftarrow{(f,\alpha)} (Y,d)$ in $\Gr(\C)$,
  consisting of a morphism $f : X \to Y$ in $sU^G_+$ and a morphism
  $\alpha : \C(\phi)c \to d$ in $\C(Y)$, to the morphism 
  $((X,c) \xleftarrow{(f,\gamma_2(f,\alpha))} (Y,d))$ in $\Gr(\ho \C)$,
  where $\gamma_2(f,\alpha)$ is the composite $\ho \C(\phi) c
  \xrightarrow t \C(\phi) c \xrightarrow \alpha d$. Here $t$ is part of
  the structure of the total left derived functor of $\C(\phi)$. 
  Using the universal property of total left derived functors it is
  easy to see that $\gamma$ preserves identity morphisms and
  compositions.   
\end{proof}
In the case where we are given a commutative semi-ring in $\C_0$ we can
combine 
Proposition \ref{monoidtam} and Construction \ref{commonsection} to
obtain a partial functor $sU^G_+ \to \Gr(\C)^{\op} \to \Gr(\ho
\C)^{\op}$. On the other hand, given a commutative semi-ring in
$\C^{\op}_0$ we would like to obtain a partial functor $\sigma : sU^G_+
\to \Gr(\ho \C^{\op})^{\op}$. This will be possible if the map $t :
\ho \C(\phi)c \to \C(\phi)c$ is an isomorphism for every object $c$ in
the image of $\sigma$ and every $\phi \in U^G_+$. Note that in this
case we can even arrange for $t :
\ho \C(\phi)c \to \C(\phi)c$ to be an identity morphism by changing
our choice of the total left derived functor $\ho \C(\phi)$ of
$\C(\phi)$. In \ref{spectrf} we consider a functor $\sigma :
sU^G_+ \to 
\Gr(\ho \C^{\op})^{\op}$ which is not constructed from a commutative semi-ring
in $\C^{\op}$. 
\section{Homotopical Tambara Categories}
\label{sechtytam}
In this section we 
show that the Tambara categories constructed from the categories of
pointed topological spaces, orthogonal 
spectra and chain complexes come with Quillen model structures
satisfying the assumptions of Corollary 
\ref{htytambara}. In particular these categories have associated
homotopy Tambara categories.

\subsection{Pointed Topological Spaces}
\label{ptspaces}
Let $\Topp$ denote the category of compactly generated pointed
topological 
spaces. This is a closed symmetric monoidal model category with a set $I$ of
generating cofibrations given by standard
inclusions of the form $S^{n-1}_+ \to D^n_+$ for $n \ge 0$. A set $J$ of
generating acyclic cofibrations is given by the maps $D^n_+ \to (D^n
\times [0,1])_+$ induced by the lower inclusion $D^n \cong D^n \times
\{0\} \subseteq D^n \times [0,1]$.

Given a $G$-set $X$ we denote by $\mathcal L_X$ the category of
set-valued functors 
$F$ on the
translation category of $X$ with the property that
for every $x \in X$ the group $X(x,x) \subseteq G$
acts transitively on $F(x)$, and we let $L_X$ denote a set of
representatives for the isomorphism classes of $\mathcal L_X$. (If $X
= \{x\}$ is a one element $G$-set, the category $\mathcal L_X$ is
isomorphic to the category of transitive $G$-sets, and we could choose
$L_X$ to be 
the family of functors corresponding to the transitive $G$-sets $X/H$,
with $H$ running over the conjugacy classes of subgroups of $G$.)
We denote by $I_X$ the set of morphisms in
$\Topp(X) = \Fun(X,\Topp)$ of the form $\id \wedge \alpha : F_+ \wedge
S^{n-1}_+ \to F_+ \wedge
 D^n_+$ for $F \in L_X$ and $\alpha \in I$.
Here $F_+ \colon X \to \Topp$ is the functor with $F_+(x)$ given by the
topological space obtained by adding a 
disjoint base point to the set $F(x)$ for $x \in X$.
We denote by $J_X$ the set of morphisms in
$\Topp(X)$ of the form $\id \wedge \alpha : F_+ \wedge D^n_+ \to F_+
\wedge (D^n \times [0,1])_+$ for $F \in L$ and $\alpha \in J$.
The following theorem is well-known. See for example \cite[Section 1.2]{DK}
or \cite[Theorem III.1.8]{MM}
\begin{thm}
\label{dkthm}
  There is a model structure on $\Topp(X)$ with $I_X$ and $J_X$ as sets of
  generating cofibrations and generating acyclic cofibrations
  respectively. In this model structure a map $f : A \to B$ is a
  fibration if and only if for every
  object $x$ 
  of $X$ and every subgroup $H$ of $X(x,x)$ the map $A(x)^H \to
  B(x)^H$ of $H$-fixpoints induced by $f_x$ is a fibration  in
  $\Topp$. It is a a weak equivalence if and only if every for every
  $x$ and $H$ as above the 
  map  
  $A(x)^H \to
  B(x)^H$ is a 
  weak
  equivalence.
\end{thm}
\begin{prop}
\label{toppisok}
  Suppose that $f : X \to Y$ is a map of finite $G$-sets. 
  The map   $(f_\wedge \overline C)(-1) \to (f_\wedge \overline C)(0)$
  is a cofibration for
  every generating
  cofibration $c : C_{-1} \to C_0$ in $\Topp(X)$ and it is an acyclic
  cofibration for
  every generating
  acyclic cofibration $c : C_{-1} \to C_0$ in $\Topp(X)$.
\end{prop}
\begin{proof}
  We only prove the statement about acyclic cofibrations.
  The map $c$ is of the form $F_+ \wedge \alpha$ for an $F$ in
  $L_X$ and 
  $\alpha : A_{-1} = D^n_+ \to A_0 =
  (D^n \times [0,1])_+$ in $J$. We have that $(f_\wedge \overline
  C)(-1) \cong 
  (f_\wedge F_+) \wedge (f_\wedge \overline A)(-1)$ 
  and that $(f_\wedge \overline C)(0) \cong
  (f_\wedge F_+) \wedge (f_\wedge \overline A)(0)$. The map
  $(f_\wedge \overline A)(-1) \to (f_\wedge \overline A)(0)$ is a
  an acyclic cofibration in $\Topp$ by the push-out-product axiom,
  and there is 
  an isomorphism of the form $(f_\wedge F_+) \cong (M_1)_+ \vee \dots
  \vee (M_n)_+$ 
  with $M_1,\dots,M_n \in L_Y$.
  The result now follows from
  the fact that the functor $(f_\wedge F_+) \wedge - : \Topp \to
  \Topp(Y)$ commutes with colimits.
\end{proof}
\begin{cor}
  There is a Tambara category $\ho \Topp$ with $\ho \Topp(X)$
  given by the homotopy category of $\Topp(X)$ and with $\ho \Topp (\phi)$
  given by a chosen total left derived functor of the functor
  $\Topp(\phi)$ for $\phi
  \in sU^G_+(X,Y)$.
\end{cor}
\begin{proof}
  This follows by combining Proposition \ref{toppisok} and
  Corollary 
  \ref{htytambara}. 
\end{proof}
\begin{remark}
  The above results also hold in the context of simplicial sets. 
  The simplicial analogue of Theorem \ref{dkthm} is treated in 
  \cite[Theorem 9.5]{DRO}.
\end{remark}
\subsection{Orthogonal Spectra}
Let $V$ be a finite-dimensional real inner product space. We denote by
$S^V$ the one-point compactification of $V$ with the added point $\infty$
as base-point.
The addition in $V$ extends to $S^V$ by 
declaring $x + y = \infty$ if either $x$ or $y$ is equal to $\infty$.
In other words $x+y$ is the image of $(x,y) \in S^V
\times S^V$ under the composition $S^V \times S^V \to S^V \wedge S^V
\cong S^{V \oplus V} \to S^V$, where the last map is
induced by addition in $V$. 
The group of linear isometries of $V$ is denoted $O(V)$.
Let $V \subseteq W$ be an inclusion of finite-dimensional real inner
product spaces and let $W-V \subseteq W$ denote the orthogonal
complement of $V$. 
Given $y \in S^W$ the pointed map $t_y \colon S^V \to
S^W$ is 
defined by $t_y(x) = x+y$ for $x \in S^V$. Given $k \in O(W)$ we
extend $k$ to a pointed map $k \colon S^W \to S^W$. 

Given an inclusion $V \subseteq W$ of finite-dimensional real inner
product spaces we let
let $\os(V,W)$ denote the subspace of $\Topp(S^V,S^W)$
consisting of maps of the form $k \circ t_y$ for $y \in S^{W-V}$ and $k
\in O(W)$.
Let $y \in S^{W-V}$ and $h \in O(V)$ and pick $h' \in O(W)$ with the
property that $h'(x) = h(x)$ for $x \in V$ and such that $h'(y) = y$.
Then $t_y \circ h = h' \circ t_y$. In particular $k \circ t_y \circ h
\circ t_x = k \circ h' \circ t_{x+y} \in \os(U,W)$ for every $k \in
O(W)$, $h \in O(V)$, $y \in S^{W-V}$ and $x \in S^{V-U}$. 
It follows that the composition $\Topp(S^V,S^W) \wedge \Topp(S^U,S^V)
\to \Topp(S^U,S^W)$ induces a composition $\os(V,W) \wedge \os(U,V) \to
\os(U,W)$ making $\os$ a category enriched over $\Topp$ with the class
of finite dimensional real inner product spaces as object class. Further
the direct sum of vector spaces induces a symmetric monoidal product
$\oplus$ on $\os$.

Note that
$\Phi \colon O(W)_+ \wedge_{O(W-V)} S^{W-V}
\to \os(V,W)$
defined by $\Phi([k,y]) =k \circ t_y$ for $k \in O(W)$ and $y \in
S^{W-V}$
is a homeomorphism. 

The following definition is taken from \cite[Example 4.4]{MMSS}.
\begin{defn}
  The category $\osp$ of orthogonal spectra is the category of
  $\Topp$-functors from $\os$ to $\Topp$. 
\end{defn}
Defining the smash-product $E \wedge F$ of orthogonal spectra $E,F \in
\osp$ as the left Kan extension of the composition $\os \times \os
\xrightarrow{E \times F} \Topp \times \Topp
\xrightarrow{\wedge} \Topp$ along 
$\oplus \colon \os \times \os \to \os$ as in
\cite[Definition 21.4.]{MMSS} and defining internal function-objects
as in \cite[Definition 21.6.]{MMSS} the category $\osp$ of orthogonal
spectra becomes a complete
and cocomplete closed symmetric monoidal category. Let us abuse
notation and denote by $\osp$ the Tambara functor with
$\osp(X) = 
\Fun(X,\osp)$ for a finite $G$-set $X$.

Given a finite $G$-set
$X$ we let $\sk(\os(X^{\op}))$ denote a set containing one representative
for each isomorphism class of objects in $\os(X^{\op}) = \Fun(X^{\op},\os)$. 
For $A$ in $\os(X^{\op})$ and $V$ in $\os$ there is a functor
$\os(A,V) : X \to \Topp$ 
because $\os$ is enriched over $\Topp$.
Thus we can consider the object
$\os(A,-)$ of $\osp(X)$.
The set
$FI_X$ consists of the morphisms $\os(A,-) \wedge \alpha$
in $\osp(X)$ 
for $A \in \sk(\os(X^{\op}))$ 
and  
$\alpha$ in $I_X$ of \ref{ptspaces}.
The set
$FJ_X$ 
consists of the morphisms $\os(A,-) \wedge \alpha$
in $\osp(X)$ 
for $A \in \sk(\os(X^{\op}))$ 
and  
$\alpha$ in $J_{X}$ of \ref{ptspaces}.

For the following result we refer to \cite[Theorem
4.2]{DRO} and \cite[Lemma 6.5]{MMSS}.
\marginpar{style}
\begin{thm}
  For every finite $G$-set $X$ there is a model structure on the
  category $\osp(X)$ with $FI_X$ as set of generating 
  cofibrations and with $FJ_X$ as set of generating acyclic
  cofibrations. A morphism $\alpha : E \to B$ is a fibration in this
  model category if and only if for every object 
  $A$ of $\os$ the map $\alpha(A) : E(A) \to B(A)$ is a fibration in
  $\Topp(X)$, and it is a 
  weak equivalence if and only the map $\alpha(A)$ is a 
  weak
  equivalence in $\Topp(X)$ for every object $A$ of $\os$.
  We call this the projective model
  structure on $\osp(X)$.
\end{thm}
We refer to fibrations, weak equivalences and acyclic
cofibrations in the projective model structure on 
$\osp$
as projective fibrations, projective weak equivalences
and projective acyclic cofibrations respectively. The cofibrations in
the projective model structure will be referred to simply as
cofibrations.
\begin{lem}
\label{projcofst}
  For every map $f:X \to Y$ of finite $G$-sets 
  the map $(f_{\wedge} \overline C)(-1) \to (f_{\wedge}
  \overline C)(0)$ is a
  cofibration in $\osp(Y)$ for 
  every generating
  cofibration 
  $c: C_{-1} \to 
  C_{0}$ in
  $\osp(X)$
  and it is a projective acyclic cofibration for every 
  projective acyclic 
  $c: C_{-1} \to 
  C_{0}$ in
  $\osp(X)$.
\end{lem}
\begin{proof}
  Using the observation $f_\wedge \os(A,-) \cong \os(f_\oplus A,-)$,
  the proof is similar to the proof of Proposition \ref{toppisok}.
\end{proof}
  Combining Lemma \ref{projcofst} and
  Corollary 
  \ref{htytambara0} we obtain: 
\begin{cor}
\label{projcor}
  For every morphism $\phi \in sU^G_+$ the functor
  $\osp(\phi)$ preserves both cofibrations and projective acyclic
  cofibrations between
  cofibrant objects. 
\end{cor}

Before defining the stable equivalences in $\osp$ we discuss a general
construction on a cocomplete symmetric monoidal category
$(\C,\diamond,n_\diamond)$. Given morphisms $a: A_{-1} \to
A_0$ and $b : B_{-1} \to B_0$ in $\C$ we have the morphism $a\Box b :
(\overline A \diamond \overline B)(-1) \to 
(\overline A \diamond \overline B)(0)$. 
It is well-known that this operation
$\Box$ is a symmetric monoidal product on the category 
$\Map \C$ of arrows in $\C$ where, for $a$ and $b$ as above, a
morphism $a \to b$ in $\Map \C$ consists of maps $\phi_i : A_i \to
B_i$ for $i = -1,0$ such that $\phi_0 a = b \phi_{-1}$
(see e.g. \cite[p. 109]{Hovey}). 

Given a morphism $A \to B$ in $\os(X^{\op})$ such that for every $x \in X$
the map $S^{A(x)} \to S^{B(x)}$ is induced by an inclusion of real
inner-product spaces we let $S^{A-B}$ denote the object of $\os(X^{\op})$
with $S^{A-B}(x) = S^{A(x)-B(x)}$ for $x \in X$. Since $X$ is a
groupoid, there is an isomorphism $X \to X^{\op}$ of categories taking
a morphism to its inverse. We consider
$\os(A,B)$ as a functor from $X$ to $\Topp$ via the composition $X \to X
\times X \to X \times X^{\op} \to \Topp$ and we let 
$\lambda_{A,B}$ denote the composition 
\begin{displaymath}
  \os(B,-) \wedge S^{B-A} \subseteq
  \os(B,-)\wedge \os(A,B) 
  \to \os(A,-).
\end{displaymath}
We choose a factorization 
\begin{displaymath}
  \os(B,-) \wedge S^{B-A} 
  \xrightarrow {k_{A,B}} 
  M \lambda_{A,B} \xrightarrow {r_{A,B}} 
  \os(A,-)
\end{displaymath}
of $\lambda_{A,B}$ by a cofibration $k_{A,B}$ and a projective
acyclic fibration $r_{A,B}$.
Let $E_X$ denote the set of morphisms of the form $k_{A,B} \Box i$ where
$k_{A,B}$ is of the above form and $i \in I_X$. 
We define $K_X = FJ_X \cup E_X$.
The following is \cite[Theorem III.5.3.]{MM}.
\begin{thm}
  For every finite $G$-set $X$ there is a model structure on $\osp(X)$
  with $FI_X$ as set of generating 
  cofibrations and with $K_X$ as set of generating acyclic
  cofibrations. 
  This is the stable model structure on $\osp(X)$.
\end{thm}
We refer to fibrations, weak equivalences and acyclic
cofibrations in the stable model structure on $\osp(X)$ as stable
fibrations, stable weak equivalences 
and stable acyclic cofibrations respectively. 
We refer to \cite[Proposition III.7.5]{MM} or \cite[Lemma 6.27]{DRO} for
the following lemma:
\begin{lem}
  The pushout-product axiom holds in $\osp(X)$ with the stable model
  structure. 
\end{lem}
The next lemma is a reformulation of \cite[Proposition V.2.3]{MM}.
\begin{lem}
  Let $f : X \to Y$ be a map of finite $G$-sets. The functor $f_\vee :
  \osp(X) \to \osp(Y)$ is a left Quillen functor with respect to the
  stable model structure.
\end{lem}
The following lemma follows from \cite[Lemma V.2.2]{MM}
\begin{lem}
  Let $f : X \to Y$ be a map of finite $G$-sets. The functor $f^* :
  \osp(Y) \to \osp(X)$ is a left Quillen functor.
\end{lem}
\begin{prop}
\label{stabisok}
  For every map $f : X \to Y$ of finite $G$-sets and every $c : C_{-1}
  \to C_0$ in $K_X$ the map $(f_\wedge \overline C)(-1) \to
  (f_\wedge 
  \overline C)(0)$ is an acyclic cofibration in the stable model
  structure on $\osp(Y)$.
\end{prop}
\begin{proof}
  If $c \in FJ_X$ we
  have by \ref{projcofst} that $f_\Box c : (f_\wedge \overline C)(-1) \to
  (f_\wedge 
  \overline C)(0)$ is a projective acyclic
  cofibration, and in particular it is a stable acyclic cofibration.
  Thus we are left with the case where $c$ is in $E_X$.
  Assume by 
  induction that the proposition holds if $f$ has fibers of
  cardinality less than or equal to $n-1$.
  Suppose that $f$ has fibers of
  cardinality less than or equal to $n$.
  By the push-out product axiom in $\osp(Y)$ we can without loss of
  generality assume that the fibers for $f$ all have cardinality  $n$. 

  Let us start by considering the case 
  $c = k_{A,B}$ for ${k_{A,B}} \colon \os(B,-) \wedge S^{B-A} 
  \to M \lambda_{A,B}$ 
  as above, where $C_{-1} = \os(B,-) \wedge S^{B-A}$
  and  $C_0 
  = M \lambda_{A,B}$. 
  With the
  notation of \ref{hatdefn} we
  have a commutative 
  diagram of the form
  \begin{displaymath}
  \begin{CD}
    (f_\wedge \widehat C)(-n) @>{\gamma}>> (f_\wedge \widehat C)(-1) 
  @>{{f_\Box k_{A,B}}}>> (f_\wedge \widehat C)(0) \\
  @V{\cong}VV @. @V{\delta}VV \\
  f_\wedge (C_{-1}) @>{\cong}>> \os(f_\oplus B,-)
  \wedge S^{f_\oplus(B-A)} @>{\lambda_{f_\oplus A, f_\oplus B}}>>
  \os(f_\oplus A,-). 
  \end{CD}
  \end{displaymath}
  Using our inductive assumption and Corollary \ref{projcor} we can
  apply Lemma 
  \ref{usepop} to conclude 
  that $\gamma$ is a stable acyclic cofibration. Since
  $M\lambda_{A,B}$ and 
  $\os(A,-)$
  are cofibrant it follows by combining Lemma \ref{projcor}
  and 
  K. Brown's lemma (see e.g. \cite[Lemma 9.9]{Dwyer-Spalinski}) 
  that $\delta$ is a projective equivalence. In particular $\delta$ is
  a stable equivalence.
  The lower horizontal map 
  $\lambda_{f_\oplus A, f_\oplus B}$
  is a stable
  equivalence by definition. It follows that ${f_\Box k_{A,B}}$ is a
  stable equivalence. Since the cofibrations
  in the projective- and in the stable model structure are the same we
  can conclude by Lemma \ref{projcor} 
  that ${f_\Box k_{A,B}}$ is a stable acyclic cofibration. 
  For
  $c$ of the form $c = k_{A,B} \Box i$ the associativity and
  commutativity isomorphisms for $\Box$ induce an isomorphism 
  $f_\Box c
  \cong (f_\Box k_{A,B}) \Box (f_\Box i)$.
  This is an acyclic cofibration by the pushout-product
  axiom. We conclude that for every $c \in E_X$ the map
  $f_\Box c$ is an acyclic cofibration. 
\end{proof}
Applying Proposition \ref{stabisok}, Corollary \ref{htytambara0},
Corollary \ref{htytambara} 
and K. Brown's lemma (see e.g. \cite[Lemma 9.9]{Dwyer-Spalinski}) 
we obtain:
\begin{cor}
  For every morphism $\phi \in sU^G_+$ the functor $\osp(\phi)$
  preserves stable equivalences between cofibrant objects. In
  particular there is a Tambara category $\ho \osp$ with $\ho
  \osp(\phi)$ given by a (chosen) total left derived functor of
  $\osp(\phi)$. 
\end{cor}
\subsection{Chain Complexes}
Let us consider the category $\Ch_R$ of chain complexes over a
commutative 
ring $R$ with the projective model structure, that is, the model
structure where the weak equivalences are the quasi-isomorphisms and
where the fibrations are the surjective chain-homomorphisms. 
Given a $G$-set $X$ we consider the model structure on the category
$\Ch_R(X) = \Fun(X,\Ch_R)$ where a map $\alpha : A \to B$ is a
fibration or weak equivalence if and only if for every $x \in X$
the map $\alpha_x : A(x) \to B(x)$ is so. It is well-known that such a
model structure exists (compare e.g. to \cite[Theorem 6.5]{MMSS}).
In particular $\Ch_R(G/G)$ is isomorphic to the category
$\Ch_{R[G]}$ of chain complexes over the group-ring $R[G]$ with the
projective model structure.
Let us consider the forgetful functor $U : \F_G \to \F$ from finite
$G$-sets to finite sets. Given a finite $G$-set $X$ there is an inclusion
$j_X: UX \subseteq X$ of translation categories. Here $UX$ is 
a category with only identity morphisms.
There is an induced functor $j^*_X : \Ch_R(X) \to \Fun(UX,\Ch_R) =
\prod_{x \in X} \Ch_R$, and a morphism $\alpha$ in $\Ch_R(X)$ is a
weak equivalence or fibration if and only if $j_X^*(\alpha)$ is so.
\begin{prop}
\label{chain61}
  For every $\phi \in sU^G_+(X,Y)$ the functor $\Ch_R(\phi)$ preserves
  weak equivalences between objects $A$ in $\Ch_R(X)$ with the property that
  $j_X^*(A)$ is cofibrant in $\Ch_R(UX)$. In particular, if $j_X^*(A)$ is
  cofibrant in $\Ch_R(UX)$, then the object $j_Y^*(\Ch_R(\phi)A)$ is
  cofibrant in 
  $\Ch_R(UY)$. 
\end{prop}
\begin{proof}
  Before we treat the general case, let us for a moment assume that
  $G$ is the trivial group. Then $U$ is simply the identity functor on
  $\F$. Given a map
  $f : X \to Y$ of finite sets, the functor $f_\otimes : \Ch_R(X) \to
  \Ch_R(Y)$ is just given by iterated applications of the tensor product
  $\otimes : \Ch_R \times \Ch_R \to \Ch_R$. It follows from the
  push-out product axiom for $\Ch_R$ that the
  functor $f_\otimes$ preserves weak equivalences between cofibrant
  objects. The functors $f_\oplus : \Ch_R(X) \to \Ch_R(Y)$ and
  $f^*:\Ch_R(Y) \to \Ch_R(X)$ are easily seen to preserve weak
  equivalences between cofibrant objects. Hence the theorem holds in
  the case where $G$ is the trivial group.

  Let us return to the general case where $G$ is a finite group. Given
  a map $f : X \to Y$ of finite $G$-sets we obtain a commutative
  diagram of the form
  \begin{displaymath}
    \begin{CD}
      \Ch_R(X) @>j^*_X>> \Ch_R(UX) \\
      @V{\Ch_R(\phi)}VV @V{\Ch_R(U\phi)}VV \\
      \Ch_R(Y) @>j^*_Y>> \Ch_R(UY).
    \end{CD}
  \end{displaymath}
  Since the result holds for the trivial group we obtain that if
  $j_X^* A$ is cofibrant, then also $j_Y^* \Ch_R(\phi)A =
  \Ch_R(U\phi)j_X^* A$ is cofibrant. Further if $j_X^*A$ and $j_X^* B$
  are cofibrant and $\alpha : A \to B$ is a weak equivalence, then
  also $j^*_X \alpha$ and $\Ch_R(U \phi)(j^*_X \alpha) = j_Y^*
  \Ch_R(\phi)(\alpha)$ are weak equivalences and can conclude that $
  \Ch_R(\phi)(\alpha)$ is a weak equivalence.
\end{proof}
Combining Proposition \ref{chain61} and Proposition \ref{derivedtambara} we
get: 
\begin{cor}
  There exists a
  Tambara category $\ho \Ch_R$ with $\ho \Ch_R(\phi)$ given by a
  total left derived functor of $\Ch_R(\phi)$ for $\phi \in sU^G_+$.
\end{cor}
\section{Homotopical Tambara Functors}
\label{sechotrf}
In this section we shall use the Tambara categories of Section
\ref{sechtytam} 
to construct Tambara functors. In order to do so we apply
the construction of Proposition
\ref{monoidtam}. 

\subsection{Transfer for Chain Complexes}
\label{evens-example}
Let us return to the situation of example
\ref{evens-example1}, where we considered the the category $\Ch_R$ of
chain complexes over a commutative ring $R$ and the Tambara category
$\Ch_R$ with $\Ch_R(X) = \Fun(X,\Ch_R)$. We saw in Proposition \ref{chain61}
that the 
Tambara category $\Ch_R$ satisfies the assumptions of Proposition
\ref{derivedtambara}, and therefore,
combining 
Proposition \ref{monoidtam} and Construction \ref{commonsection},
a differential graded commutative
$R$-algebra 
$A$ gives rise to a functor $\sigma_A$ given by the composition
$sU^G_+ \to \Gr(\Ch_R)^\op \to \Gr(\ho \Ch_R)^\op$. 
On the other hand, 
if $B$ is a differential graded cocommutative $R$-coalgebra and if $B$
is cofibrant considered as a chain complex of $R$-modules, then by
Proposition \ref{chain61}, for
every $\phi \in U^G_+(X,Y)$ the map $t : \ho \Ch_R(\phi)p_X^* B \to
\Ch_R(\phi) p_X^* B$ is an isomorphism in $\ho \Ch_R(Y)$, and by
Construction \ref{commonsection} we obtain a partial functor $\sigma_B :
sU^G_+ \to 
\Gr(\ho \Ch_R^\op)^\op$ with $\sigma_B(X) = (X,p_X^*B)$ and with
$\sigma_B(\phi:X \to Y)=(\phi,\sigma_{B,2}(\phi))$ where $\sigma_{B,2}$
is given by the composite $\ho \Ch_R(\phi) p_X^*B \xleftarrow {t^{-1}}
\Ch_R(\phi) p_X^* B \xleftarrow{s_{B,2}(\phi)} p_Y^* B$ for $s_B(\phi) =
(\phi,s_{B,2}(\phi))$ of \ref{commonsection}. Combining the functors 
$\sigma_A$ and $\sigma_B$ we obtain a product-preserving partial
functor $sU^G_+ \to 
\Gr(\ho \Ch_R^\op)^\op \times_{sU^G_+} \Gr(\ho \Ch_R)^\op \to \Set$,
and by \ref{simpleenough} this partial functor extends to a Tambara
functor $\ho\Ch_R(B,A) : U^G_+ \to \Set$. In the particular case where
$B = \Z[T]$ with $T$ in degree $2$ and $A = R$ we obtain a Tambara
functor $\ho \Ch_R(\Z[T],R)$ with $\ho \Ch_R(\Z[T],R)(G/H) \cong
H^*(BH,R) = H^*(H,R)$. This Tambara functor was considered by Tambara
in \cite[3.4]{Tambara}.

\subsection{Transfer for Orthogonal Spectra}
\label{spectrf}
We associate a Tambara functor
$\widetilde A$
with
$$\widetilde A(X) = 
\ho
\osp(X)(p_X^*\So,p_X^*A) \cong
\ho
\osp(*)((p_X)_{\vee} p_X^* \So,A)
=
[\Sigma^{\infty}X_+,A]_G
$$
to every commutative
orthogonal $G$-ring\-spectrum $A$, that is,
to every commutative monoid $A$ in $\osp(*)$. 
Here 
$\So \in \osp(*)$ is the
$G$-sphere spectrum defined by $\So(V) = S^V$ for $V \in \os(*)$ and the map
$p_X : X \to *$ is the unique map to the terminal object in
$\FG$. We use the notation $(p_X)_{\vee} = \osp(T_{p_X})$ and $p_X^* =
\osp(R_{p_X})$.  
Note that
for $X = G/H$ we have  
$
\widetilde A (G/H) = [\Sigma^{\infty} G/H_+,A]_G = \pi_0(A^H),
$
where $A^H$ denotes the $H$-fixed point spectrum of $A$.

The construction of $\widetilde A$ is complicated by the fact
that the sphere spectrum $\So$ is a coalgebra  in $\ho \osp$ but 
not in $\osp$. Note that we do not require $A$ to be cofibrant.

In \ref{commonsection} we constructed a section 
$s_A : 
sU^{G}_+ \to \Gr(\osp)^{\op}$ of the partial functor
$\GR(\osp)^{\op} : \Gr(\osp)^{\op} \to sU^{G}_+$. Composing with the
functor $\Gr(\osp)^{\op} \to \Gr(\ho\osp)^{\op}$ we obtain a
product-preserving section
$\sigma_A : sU^{G}_+ \to \Gr(\ho \osp)^{\op}$ of the
partial functor
$\GR(\ho \osp)^{\op} : \Gr(\ho \osp)^{\op} \to sU^{G}_+$.
In this section we construct a product-preserving section 
$\sigma_{\So} \colon sU^G_+
  \to \Gr((\ho \osp)^{\op})^\op$ of the functor
$\GR((\ho \osp)^{\op})^{\op}$.
 
The composite 
\begin{displaymath}
  sU^{G}_+ \xrightarrow {(\sigma_{\So},\sigma_A)}  (\Gr((\ho \osp)^{\op})
  \times_{(sU^{G}_+)^{\op}} \Gr(\ho \osp))^{\op} 
  \xrightarrow {\hom_{\ho \osp}} \Set.
\end{displaymath}
is a product-preserving partial functor, and by \ref{simpleenough} 
it defines a Tambara functor $\widetilde A = \ho \osp(\So,A)$
with
\begin{eqnarray*}
  \widetilde A (X) = \ho \osp(X)(p_X^*\So,p_X^*A) 
  \cong [\Sigma^{\infty} X_+,A]_G
\end{eqnarray*}
for $X$ in $\FG$. 
The rest of this section contains a construction of $\sigma_{\So}$.

Given a symmetric monoidal category $(\C,\diamond,n_{\diamond})$
and $f: X \to Y$ in $\F_G$ the functor $f_\diamond \colon \Fun(X,\C)
\to \Fun(Y,\C)$ with $f_\diamond(c)(y) \cong
\Diamond_{x \in f^{-1}(y)} c(x)$ for $c \in \Fun(X,\C)$ is constructed
similarly to $\C(N_f)$ of Construction \ref{theTamcat}.
In particular we use this notation for $\C = \Topp$ with the monoidal
products $\amalg$, $\times$ and $\wedge$.

For typographical reasons we introduce the notation $\D = (\ho \osp(X))^\op$.
  A functor $\sigma = (\id,\sigma_2) : sU^G_+
  \to \Gr(\D)^\op$ is uniquely
  determined by its values $\sigma(R_f)$, $\sigma(T_f)$ and
  $\sigma(N_f)$ for $f \colon X \to Y$ in $\F_G$, and given such values,
  they extend to a functor $\sigma$ if and
  only if the generating relations between $R_f,T_f$ and $N_f$ in
  $sU^G_+$ considered in \cite[Proposition 7.2]{Tambara} are
  respected. Below we define such a functor $\sigma$ by specifying these values
  and verifying that the
  relations of \cite{Tambara} are satisfied.
Since it is more convenient we work with a strictly unital
smash-product in $\osp$. In particular this implies
$\So_X = \osp(N_{i_X})(*)$ where $i_X :
\emptyset \to X$ is the unique map from the initial object of $\FG$ to
$X$.   
For $\phi
\in sU^G_+(X,Y)$ we choose the total left derived functor $\ho
\osp (\phi)$ such that if $M$ is a cofibrant object in $\osp(X)$,
then $t : \ho 
\osp (\phi)M \to \osp(\phi)M$ is the identity map. 

  Given $f : X \to Y$ in $\F_G$ we have
  \begin{displaymath}
    \D(R_f) \So_Y = \D(R_f) \D(R_{p_Y}) \So
    = \D(R_f R_{p_Y}) \So = \D(R_{p_X})\So = \So_X, 
  \end{displaymath}
  and we let $\sigma(R_f) = (R_f ,\id)$. Similarly we have
  \begin{displaymath}
    \D(N_f) \So_X = \D(N_f) \D(N_{i_X})*
    = \D(N_f N_{i_X}) * = \D(N_{i_Y})*  = \So_Y.
  \end{displaymath}
  and we let $\sigma(N_f) = (N_f,\id)$.
  We let $\sigma(T_f) = (T_f,t_f)$, where the map $t_f$ is the
  transfer map defined below.
  If $f$
  is the projection 
  $G/H \to G/G$, then $t_f$ is the classical transfer map
  of Kahn and Priddy
  \cite{Kahn-Priddy} and of Roush \cite{Roush}. This map has been studied
  further for example by Becker and Gottlieb
  \cite{Becker-Gottlieb}. In order to 
  give the definition of the 
  transfer map in our context we note that the fibers
  $f^{-1}(y)$ for $y \in Y$ assemble to a functor $f^{-1} : Y \to
  \Ens$. We can choose a functor $V : Y \to \os$ and an embedding
  $\iota_f \colon f^{-1} \hookrightarrow V$, that is, embeddings
  $\iota_{f,y} \colon f^{-1}(y) 
  \hookrightarrow V(y)$ for $y \in Y$ such that
  $\iota_{f,{gy}}\circ f^{-1}(g) = V(g) \circ \iota_{f,y}$
  for every $y \in Y$ and $g \in G$.
  Identifying small disjoint balls around the
  elements of $f^{-1}(y)$ with $V(y)$
  the embedding $\iota_f$ can be extended to an embedding $\iota_f' \colon
  f^{-1} \times V \hookrightarrow V$. Collapsing the complement of 
  $f^{-1}(y) \times V(y)$ in $V(y)$ for $y \in
  Y$ we obtain a map 
  $$\tau_{f,y} \colon S^{V(y)} \to
  (f_\vee f^* S^V)(y) = \bigvee_{x \in f^{-1}(y)} S^{V(y)}.$$
  These maps define a map
  $\tau_f
  : S^V \to 
  f_{\vee} f^* S^V$.
  We denote the suspension spectrum of $Z \in \Fun(Y,\Topp)$
  by $\So_Y \wedge Z$.
  We define the map $t_f \colon \So_Y \to f_\vee \So_X = \D(T_f)\So_X$
  in $\ho \osp(Y)$ as the
  composition:
  \begin{eqnarray*}
    \So_Y &\xrightarrow \cong& \hom(\So_Y \wedge S^V, \So_Y \wedge
    S^V) \\
    &\xrightarrow {(\tau_f)_*}& \hom(\So_Y \wedge S^V, \So_Y \wedge
    f_\vee f^*S^V) \\
    &\xrightarrow \cong& \hom(\So_Y \wedge S^V, f_\vee f^* \So_Y \wedge
    S^V) \\
    &\xrightarrow \cong& f_\vee f^* \So_Y = f_\vee \So_X.
  \end{eqnarray*}
  Here $\hom$ denotes
  the internal hom-object in $\ho \osp(Y)$.
  In other words we have
  defined a map $t_f \colon \D(T_f)S_X \to S_Y$ in $\D(Y)$.
  This ends the construction of the morphisms
  $\sigma(R_f)$, $\sigma(T_f)$ and $\sigma(N_f)$ for $f$ a morphism in
  $\FG$. 
  We move on to the verification of the relations between 
  them.

  Note first that $t_f = t_f(\iota_f)$ is
  independent of the embedding $\iota_f \colon f^{-1}
  \hookrightarrow V$.
  One way to see this is to note that if $\iota_1 \colon f^{-1}
  \hookrightarrow V_1$ and $\iota_2 \colon f^{-1} \hookrightarrow V_2$
  are 
  embeddings, we obtain a diagonal 
  embedding $\iota_1 \times \iota_2$ by the composition $f^{-1}
  \hookrightarrow f^{-1} \times f^{-1} 
  \hookrightarrow V_1 
  \oplus V_2$, and the diagram 
  \begin{displaymath}
    \begin{CD}
      S^{V_1 \oplus V_2} @>{\cong}>> S^{V_1} \wedge S^{V_2} \\
      @V{t_f(\iota_1 \times \iota_2)}VV @V{\id \wedge t_f(\iota_2)}VV \\
      f_\vee f^* S^{V_1 \oplus V_2} @>{\cong}>> S^{V_1}
      \wedge  (f_\vee f^* S^{V_2})
    \end{CD}
  \end{displaymath}
  is homotopy-commutative by a homotopy that shrinks the
  $V_1$-coordinate. 
  
  Since
  the relations of \cite[Proposition 7.2]{Tambara} not
  involving morphisms of the form $T_f$ are readily
  verified we concentrate on the ones involving $T_f$. 
  We first show that given a pull-back diagram of the form
  \begin{displaymath}
    \begin{CD}
      X @<g'<< X' \\
      @VfVV @Vf'VV \\
      Y @<g<< Y'
    \end{CD}
  \end{displaymath} 
  in $\FG$ the relation 
  $\sigma(R_g) \sigma(T_f) = \sigma(T_{f'})
  \sigma(R_{g'})$ holds. This translates into showing that 
  the diagram 
  \begin{displaymath}
    \begin{CD}
      \D(R_g) \D(T_f) \So_X @>{\D(R_g)(t_f)}>> \D(R_g) \So_Y
      @>{\id}>> \So_{Y'} \\
      @V{\cong}VV @. @| \\
      \D(T_{f'}) \D(R_{g'}) \So_X @>{\id}>>
      \D(T_{f'}) \So_{X'}  
      @>{t_{f'}}>> \So_{Y'}
    \end{CD}
  \end{displaymath}
  commutes.
  In order to do so we choose an embedding $f^{-1} \hookrightarrow
  V$ and consider the induced embedding ${f'}^{-1} \cong g^* 
  f^{-1} \hookrightarrow g^* V$. Tracing back the definition of the
  transfer using these embeddings we see that it suffices to
  note that the following diagram in $\Topp$ is commutative:
  \begin{displaymath}
    \begin{CD}
      g^* f_{\vee} f^* S^V @<{g^*(\tau_f)}<< g^* S^V
      @= S^{g^*V} \\
      @A{\cong}AA @. @| \\
      {f'}_\vee {g'}^* f^* S^V @= {f'}_\vee {f'}^* S^{g^* V} 
      @<{\tau_{f'}}<<
      S^{g^*V} .
    \end{CD}
  \end{displaymath}

  Next we show that $\sigma(T_h) \sigma(T_f) = \sigma(T_{hf})$ for $f
  \colon X \to Y$ and $h \colon Y \to Z$ in $\FG$. This amounts to
  showing that the diagram
  \begin{displaymath}
    \begin{CD}
      \D(T_h) \D(T_f) \So_X @>{\D(T_h)(t_f)}>> \D(T_h) \So_{Y} \\
      @V{\cong}VV @V{t_h}VV \\
      \D(T_{hf}) \So_X  @>{t_{hf}}>>  \So_Z 
    \end{CD}
  \end{displaymath}
  commutes. In order to do so we choose $V : Z \to \os$ 
  and embeddings $\iota_{hf} \colon ({hf})^{-1} \hookrightarrow V$ and
  $\iota_h \colon h^{-1}
  \hookrightarrow V$. We let $\iota_f$ denote the embedding
  $f^{-1} \hookrightarrow h^* ({hf})^{-1} \hookrightarrow h^*
  V$. Tracing back the definition of the transfer with respect to these
  embeddings we see that it suffices to note that the diagram
  \begin{displaymath}
    \begin{CD}
      h_{\vee} f_{\vee} f^* h^* S^V @<{h_\vee(\tau_f)}<< h_\vee
      h^* S^V \\
      @A{\cong}AA @A{\tau_h}AA \\
      (hf)_\vee (hf)^* S^V  @<{\tau_{hf}}<<  S^V
    \end{CD}
  \end{displaymath}
  commutes.

  Finally we verify the relation $\sigma(T_q) \sigma(N_{f'})
  \sigma(R_e) = \sigma(N_f) \sigma(T_p)$ for the exponential diagram of
  Definition \ref{defexpdiag}.
  This amounts to showing that the
  following diagram commutes:
  \begin{displaymath}
    \begin{CD}
      \D(T_q) \D(N_{f'}) \D(R_e) \So_A @>{\cong}>> \D(N_f) \D(T_p)
      \So_A @>{\D(N_f)(t_p)}>> \D(N_f) \So_X \\
      @VV{\id}V @.
      @V{\id}VV \\
      \D(T_q)\D(N_{f'})\So_{X'} @>{\id}>>
      \D(T_q)\So_{Y'} @>{t_q}>> \So_Y. \\ 
    \end{CD}
  \end{displaymath}
  In order to do this we choose an embedding $\iota_p \colon p^{-1}
  \hookrightarrow  V$ and we
  let $W = f_\oplus V$.
  Let $\iota_q$ denote the
  induced embedding $q^{-1} 
  \cong f_{\times} p^{-1} \subseteq f_{\times} V = f_{\oplus} V =
  W$. Tracing back definitions it suffices to note that we have a
  commutative diagram of the form:
  \begin{displaymath}
    \begin{CD}
      q_\vee {f'}_\wedge e^* (p^* S^V)  @<{\cong}<< f_\wedge (p_\vee
      p^* S^V) @<{f_\wedge(\tau_p)}<< f_\wedge S^V  \\
      @A{\cong}AA  @.
      @| \\
      q_{\vee} q^* f_\wedge S^V
      @=
      q_\vee q^* S^W @<{\tau_q}<< S^W. \\ 
    \end{CD}
  \end{displaymath}

\section{Equivariant cobordism}
\label{witam}
As an application of the our theory we consider
equivariant cobordism and the related equivariant spectrum $MU$. 
We start by recollecting 
the following theorem which was proved by Strickland in the context of
non-equivariant spectra \cite{Strickland}.
\begin{thm}
\label{MP}
  There is a commutative orthogonal ring spectrum $MP$ with action of
  $G$ of the 
  homotopy type $\bigvee _{r \in \Z} \Sigma^{2r} MU$ considered as
  a monoid in the $G$-equivariant stable category.
\end{thm}
Actually we will construct a unitary ring spectrum rather
than an orthogonal ring spectrum. There is
strong symmetric monoidal functor from the category unitary spectra to
the category of orthogonal spectra defined by left Kan extension 
\cite[Proposition 3.3]{MMSS}.
Therefore we can pass
from commutative unitary ring spectra to commutative orthogonal ring
spectra. We restrict our attention to unitary spectra
in this section.

Given a complex  inner product space $V$ and $n \in \N$ we let $\Grass(n,V)$
denote the Grassmannian manifold of (complex) $n$-dimensional subspaces of
$V$. This defines a functor $\Grass(n,-)$ from the category of inner
product spaces and injective homomorphisms to the category of
topological spaces.
Given an inclusion $V \subseteq W$, the induced map $\Grass(n,V) \to
\Grass(n,W)$ is $(2\dim(V)-c)$-connected for a constant $c$ not
depending on $V$ and $W$.
We denote by $E(n,V)$ the tautological $n$-plane bundle over
$\Grass(n,V)$ consisting of pairs $(X,x)$ of an $n$-dimensional
subspace $X$ of $V$ and a point $x \in X$. 
The associated Thom space is denoted $T(n,V)$.
If $G$ acts on $V$, then there is an induced action of $G$ on $T(n,V)$. 
For complex inner product spaces $V$ and $W$ there is a pairing 
$T(n,V)
\wedge T(m,W) \to T(n+m,V \oplus W)$
induced by the pairing
$E(n,V)
\oplus E(m,W) \to E(n+m,V \oplus W)$
taking $((X,x),(Y,y))$ to $(X \oplus Y, x+y)$.

Let $\U$ be a complete complex $G$-universe. Given an inner product space
$V$ with action of $G$ we follow tom Dieck \cite{Dieck-topology} (see
also \cite{Constenoble}) and let $MU(V) = T(\dim(V),V \oplus \U)$.
The pairing 
\begin{eqnarray*}
T(\dim(V), V \oplus \U) \wedge T(\dim(W),W \oplus \U) &\to& 
T(\dim(V \oplus W), V \oplus \U \oplus W \oplus \U) \\
&\cong& 
T(\dim(V \oplus W), V \oplus W \oplus \U)  
\end{eqnarray*}
defines a map from $MU(V) \wedge MU(W)$ to $MU(V \oplus W)$ and the
inclusion $V \to E(\dim(V),V \oplus \U)$ 
taking $v$ to $(V \oplus 0,v)$
defines a map $S^V \to MU(V)$. This structure
defines a commutative ring $G$-prespectrum $MU$ in the sense
that the prespectrum $MU$ represents a 
commutative monoid in the $G$-equivariant stable category.
For $r \ge 0$ we let $MU_r(V) = MU(V \oplus \CC^r)$ and
$MU_{-r}(V) = \map(S^{\CC^r},MU(V))$. Thus for $r \in \Z$ there is a stable
equivalences $MU_r \simeq \Sigma^{2r} MU$ of $G$-prespectra. 

\begin{proof}[Proof of Theorem \ref{MP}]
Given a finite-dimensional inner product space $V$ 
we follow Strickland \cite{Strickland} and 
define $MP(V) = \bigvee_{r \in \Z} MP_r(V)$ where 
$MP_r(V)$ denotes the space $T(r+\dim(V),V \oplus V)$ for $r \in \Z$. 
The map 
\begin{eqnarray*}
&& T(r+\dim(V), V \oplus V) \wedge T(s+\dim(W),W \oplus W) \to \\
&& T(r+s+\dim(V \oplus W), V \oplus V \oplus W \oplus W) 
\end{eqnarray*}
defines a pairing 
$MP_r(V) \wedge MP_s(W) \to MP_{r+s}(V \oplus W)$
and the
map from $V$ to $E(\dim(V),V \oplus V)$ 
taking $v$ to $(V \oplus 0,v)$
defines a map $S^V \to MP_0(V)$. 
Together these maps make $MP$ into a commutative unitary ring
spectrum. (See e.g. \cite[Theorem
  3.4]{MM} for the translation between unitary ring spectra and this
kind of 
structure.)
If $G$ acts on the inner product space $V$ then we obtain an induced
action on $MP(V)$, and this way $MP$ is considered as a commutative
unitary ring spectrum with action of $G$.

We need to show that $MP$ represents $\bigvee_{r \in \Z} 
\Sigma^{2r} \wedge MU$ as a monoid in the $G$-equivariant stable
category.
If $r \ge 0$ and $V$ is
contained in 
a $G$-universe $\U$, then the composition
\begin{eqnarray*}
  T(r+\dim(V), V \oplus V) &\to& T(r+\dim(V), V \oplus \U)
  \\
  &\to& T(\dim(V \oplus \CC^r), V \oplus \CC^r \oplus \U) 
\end{eqnarray*}
defines a map
$MP_r(V) \to MU(V \oplus \CC^r)$ and the composition
\begin{eqnarray*}
  S^{\CC^{r}} \wedge T(-r+\dim(V),V \oplus V) &\to&
  T(r,\CC^r\oplus \U) \wedge T(-r+\dim(V),V \oplus \U) \\
  &\to&
  T(\dim(V),\CC^r \oplus V \oplus \U) \\
  &\to & T(\dim(V), V \oplus \U)
\end{eqnarray*}
is adjoint to a map $MP_{-r}(V) \to \map(S^{\CC^{r}},MU(V))$.
The map 
$MP_r(V)^H \to MU_r(V)^H$  
is $(2\dim(V^H) - c)$-connected for a constant $c$ independent of $V$
for every $r \in \Z$ and every subgroup $H$ of $G$. In
particular the map 
$MP_r \to MU_r$ of prespectra
is a $\pi_*$-equivalence and it induces an isomorphism in the
$G$-equivariant  stable category.
We leave it to the reader to check that these maps are compatible
with the ring structures.
\end{proof}

In the rest of this section we prove Theorem
\ref{indthmfirst} and Theorem \ref{indthmlast}.
Let $A$ be a commutative ring and let $\CO$ denote the partially
ordered set of
conjugacy classes of subgroups of $G$ with $[H] \le [K]$ if and only there
exists $g \in G$ such that $H \subseteq gKg^{-1}$. The ring $\W_G(A)$ of
\cite{Dress-Siebeneicher} has $\map(\CO,A)$ as underlying set  and 
ring-structure 
defined through the ghost-coordinates
$\Phi_{[H]}^A 
: \W_G(A) \to A$ 
with
\begin{displaymath}
  \Phi_{[H]}^A(\alpha) = \sum_{[K] \in \CO} |(G/K)^H| \cdot 
  \alpha([K])^{(K:H)} 
\end{displaymath}
for $[H] \in \CO$ and $\alpha \colon \CO \to A$.
Here $|(G/K)^H|$ denotes
the cardinality of the set of $H$-fixed points of $G/K$
and $(K:H)$ denotes the index of $H$ in $gKg^{-1}$. 
These ghost coordinates assemble to a ring homomorphism $\Phi^A \colon
\W_G(A) \to \map(\CO,A)$ called the ghost map. 
If the underlying abelian group of $A$
is torsion free then the ghost map is injective 

In Section \ref{spectrf} we constructed the Tambara functor
$\widetilde {MP}$.
The homomorphism
$\tau_{\widetilde {MP}} : \W_G([\So,MP]) \to \widetilde {MP} (*) \cong
[\So,MP]_G$
defined by the formula
$$\tau_{\widetilde {MP}} = \sum_{[K] \in \CO}
\widetilde {MP}(T_{p^G_K} \circ N_{p^K_e}) $$
is called the {\em Teichm\"uller homomorphism} in \cite{Brun}. Here
$p^G_K$ and $p^K_e$ denote the $G$-maps  
$G/e \xrightarrow{p^K_e} G/K \xrightarrow{p^G_K} G/G$.
\begin{thm}
\label{factorghost}
  There exists a homomorphism 
  $$R \colon [\So,MP]_G \to
  \map(\CO,[\So,MP])$$ 
  of commutative rings
  with $R \circ \tau_{\widetilde {MP}} = \Phi^{[\So,MP]}$.
\end{thm}
\begin{proof}[Proof of Theorem \ref{indthmfirst}]
  By Milnor and Novikov's calculation
  \cite{Milnor},\cite{Novikov} the underlying abelian group of
  $[\So,MP] = \pi_0(MP)$ is torsion free. It follows that the ghost map
  $\Phi^{[\So,MP]}$ 
  is injective and hence by Theorem \ref{factorghost}
  $\tau_{\widetilde {MP}}$ is injective.
\end{proof}
In order to construct the homomorphism $R$ it is convenient to note
that $MP(V)$ is the Thom space of the 
canonical bundle over the space $\coprod_{r \in
\Z} \Grass(r, V \oplus V)$ of subspaces of $V \oplus V$.
Let $H$ be a subgroup of $G$ 
and $\U$ denote a complete
$G$-universe. 
An element of $MP(V)^H$ consists of a
pair $(X,x)$ of an
$H$-invariant subspace $X$ of $V \oplus V$ and an element $x \in
X^H$. To this pair we associate the pair $(X^H,x)$ considered as an
element of $MP(V^H)$. 
This defines a map $r^H \colon MP(V)^H \to
MP(V^H)$. 
As defined in \cite[Definition III.3.2]{MM} the group $[\So,MP]_H$ is
the colimit over
all finite 
dimensional subspaces $V \subseteq \U$
of the abelian groups $\pi_0(\map(S^V,MP(V))^H)$.
The composition
$$\map(S^V,MP(V))^H \to \map(S^{V^H},MP(V)^H)
\xrightarrow{\map(S^{V^H},r^H)} 
\map(S^{V^H},MP(V^H))$$
induces a map $r^H \colon [\So,MP]_H \to [\So,MP]$. 
\begin{remark}
  In other words we have constructed a map of orthogonal spectra
  $r_G \colon \Phi^G MP \to MP$ from the geometric fixed point
  spectrum of $MP$ 
  back to $MP$. (See \cite[Section V.4]{MM} for a definition of the
  geometric fixed point spectrum.) Given $W \subseteq \U$, there is a
  map $s_G(W) \colon MP(W^G) \to MP(W)^G$ taking a pair $(x,X)$ with $X
  \subseteq W^G \oplus W^G$ to the same pair with $X$ considered as a
  subspace of $W \oplus W$. These maps define a map $s_G \colon MP \to
  \Phi^G MP$ with $r_G \circ s_G = \id$.
\end{remark}
\begin{lem}
\label{rht}
  $r^{H} \circ T_{p^H_K} = 0$ for $K < H$.
\end{lem}
\begin{proof}
  This follows from the fact that $(H/K_+ \wedge S^V)^{H} = *$.
\end{proof}
\begin{lem}
\label{rhn}
  $r^{H} \circ N_{p^H_e} = \id$ for $H \le G$.
\end{lem}
\begin{proof}
  This follows from the facts that the $H$-fixed points space of
  $(S^V)^{\wedge H}$ is $S^V$ and that the $H$-fixed point
  sub-vector space of 
  $\C[H] \otimes V = \oplus_{h \in H} V$ is isomorphic to $V$.
\end{proof}
\begin{proof}[Proof of Theorem \ref{factorghost}]
  We define $R$ by $R(a)([H]) = r^H(
  \widetilde {MP}(R_{p^G_H})(a))$ for $[H] \in \CO$ and $a \in
  [\So,MP]_G$. 
  For $\alpha \colon \CO \to [\So,MP]$ we have
  \begin{eqnarray*}
  && r^H (
  \widetilde {MP}(R_{p^G_H}) ( 
  \tau_{\widetilde {MP}}(\alpha))) 
  =
  \sum_{[K] \in \CO}
  r^H (\widetilde {MP} (R_{p^G_H} \circ T_{p^G_K} \circ
  N_{p^K_e}) (\alpha([K]))).
  \end{eqnarray*}
  Using Lemma \ref{rht} and the additive double coside formula and Lemma
  \ref{rhn} and the multiplicative double
  coside formula respectively we compute for $H \le K$ that
  \begin{eqnarray*}
  r^H (\widetilde {MP}(R_{p^G_H} \circ T_{p^G_K} \circ
  N_{p^K_e})(a))
  &=&
  |(G/K)^H|  \cdot  
  r^H (\widetilde {MP} (R_{p^K_H} \circ N_{p^K_e})(a) \\
  &=& 
  |(G/K)^H| \cdot a^{(K:H)}
  \end{eqnarray*}
  for $a \in [\So,MP]$.
  For general subgroups $H$ and $K$ of $G$ we have
  \begin{eqnarray*}
  r^H (\widetilde {MP}(R_{p^G_H} \circ T_{p^G_K} \circ
  N_{p^K_e})(a))
  = 
  |(G/K)^H| \cdot a^{(K:H)}.
  \end{eqnarray*}
  In particular
  $R \circ \tau_{\widetilde {MP}} = \Phi^{[\So,MP]}$.
\end{proof}

\begin{proof}[Proof of Theorem \ref{indthmlast}]
This proof is similar to the proof of
Theorem \ref{indthmfirst}. However we are not able to prove that
there is a Tambara functor with value $\U_*^H$ on $G/H$.
Instead
we  consider the axillary
semi-Tambara functor $\M_*$ with $\M_*(X)$ given by the set of 
isomorphism classes of functors from the translation category of
$X$ to the category of almost complex manifolds.
We let $\widehat \M_*$ denote the Tambara functor obtained by
group-completing $\M_*$. By \cite[Theorem 3.6.]{Brun} we have the
Teichm\"uller 
homomorphism $\tau_{\widehat \M_*} : \W_G(\widehat \M_*(G/e)) \to
\widehat \M_*(G/G)$. 
Choosing a ring-homomorphism $\sigma : \U_* \to \widehat 
\M_*(G/e)$ such that $\sigma$ is a section of the surjection $q:
\widehat \M_*(G/e) \to \U_*$, and 
denoting by $q^G$ the surjection $\widehat \M_*^G(G/G) \to \U^G_*$  we
obtain 
a commutative diagram of the form:
\begin{displaymath}
\begin{CD}
  \W_G(\widehat \M_*(G/e)) @>{\tau_{\widehat \M_*}}>> \widehat
  \M_*(G/G) @>{\mathrm {Fix}}_{\widehat \M_*}>> \map(\CO,
  \widehat \M_*(G/e)) \\
  @A{\W_G(\sigma)}AA @V{q^G}VV @V{\map(\CO, q)}VV \\
  \W_G(\U_*) @>{q^G \circ \tau_{\widehat \M_*} \circ \W_G(\sigma)}>>
  \U^G_* @>{\mathrm {Fix}}_{\U_*}>> \map(\CO, 
  \U_*), 
\end{CD}
\end{displaymath}
where $({\mathrm {Fix}}_{\U_*}[M])([K]) = [M^K]$ and $({\mathrm
  {Fix}}_{\widehat \M_*}([M]))([K]) = [M^K]$.
Note that we have  $\Phi^{\M_*(G/e)}
= {\mathrm {Fix}}_{\widehat \M_*} \circ \tau_{\widehat \M_*}$. Since 
$\Phi_{[K]}$ is given by an integral polynomial we have 
$$q \circ
  \Phi^{\widehat \M_*(G/e)}_{[K]} 
\circ \W_G(\sigma) =
q \circ
  \Phi^{\widehat \M_*(G/e)}_{[K]} 
\circ \map(\CO, \sigma) = \Phi_{[K]}^{\U_*} \circ \map(\CO,
q \circ \sigma) = \Phi_{[K]}^{\U_*},$$
 and therefore 
 \begin{eqnarray*}
 {\mathrm 
  {Fix}}_{\U_*} 
\circ q^G \circ \tau_{\widehat \M_*} \circ \W_G(\sigma) &=& 
\map(\CO,q) \circ {\mathrm {Fix}}_{\widehat \M_*} \circ \tau_{\widehat \M_*}
\circ \W_G(\sigma) \\
&=& \map(\CO,q) \circ \Phi^{\widehat \M_*(G/e)} 
\circ \W_G(\sigma) \\
&=& \Phi^{\U_*}.  
 \end{eqnarray*}
It now follows from the injectivity of $\Phi^{\U_*}$ 
that $\tau_{\U_*} := q^G \circ \tau_{\widehat
  \M_*} \circ \W_G(\sigma)$ is 
injective. 
\end{proof}

\section{Filtered Objects}
\label{secfilob}
We shall let $\Z$ denote the usual partially ordered set of integers
considered as a symmetric monoidal category with monoidal operation
given by the sum in $\Z$.
\begin{defn}
  A functor $X : \Z \to \C$ is called a {\em filtered object} in
  $\C$. The category  
  $\Z \C$ of functors from $\Z$ to $\C$ is the {\em category of
  filtered objects in 
  $\C$}. If $(\C,\diamond,u)$ is a cocomplete symmetric monoidal category,    
  there is a symmetric monoidal structure on $\Z \C$ induced
  from the symmetric monoidal structures on $\Z$ and $\C$ by the usual
  left Kan extension: Given $A,B : \Z \to \C$, the monoidal product
  $A\diamond B : \Z \to \C$ is the left Kan extension of the composition
  \begin{displaymath}
    \Z \times \Z \overset{A \times B} \longrightarrow \C \times \C
  \overset{\diamond} \to \C
  \end{displaymath}
  along the functor $+ : \Z \times \Z \to \Z$ with 
  $(A\diamond B)(i) = \colim{\alpha+\beta \le i} A(\alpha) \diamond B(\beta)$. 
\end{defn}
\begin{lem}
\label{cralem2b}
  Let $E \to F$ be a map of filtered objects in a symmetric
  monoidal 
  category $(\C,\diamond,u)$ satisfying that for every $c \in
\C$, the functor $c \diamond -$ preserves push-outs.
Suppose that $D$ is a filtered object in $\C$ and that for
every $i \in \Z$ the square 
  \begin{displaymath}
  \begin{CD}
    E(i) @>>> E({i+1}) \\
    @VVV @VVV \\ 
    F(i) @>>> F({i+1}) 
  \end{CD}
  \end{displaymath}
is a push-out in $\C$. Then the diagram
  \begin{displaymath}
  \begin{CD}
    (E \diamond D)(i) @>>> (E \diamond D)({i+1}) \\
    @VVV @VVV \\ 
    (F \diamond D)(i) @>>> (F \diamond D)({i+1}) 
  \end{CD}
  \end{displaymath}
  is a push-out in $\C$ for every $i \in \Z$.
\end{lem}
\begin{proof}
  We can factor the second diagram of the lemma as
  \begin{displaymath}
  \begin{CD}
    (E \diamond D)(i) @>>> 
    (F \diamond D)(i) \\
    @V{=}VV @V{=}VV \\
    \colim{\alpha_1 + \alpha_2 \le i} E(\alpha_1) \diamond D(\alpha_2)
    @>>>
    \colim{\alpha_1 + \alpha_2 \le i} F({\alpha_1}) \diamond D({\alpha_2}) \\
    @VVV @VVV \\
    \colim{\alpha_1 + \alpha_2 \le i} E({\alpha_1+1}) \diamond D({\alpha_2}) 
    @>>>
    \colim{\alpha_1 + \alpha_2 \le i} F({\alpha_1+1}) \diamond D({\alpha_2}) \\
    @V{\cong}VV @V{\cong}VV \\
    (E \diamond D)({i+1}) @>>>
    (F \diamond D)({i+1}).
  \end{CD}
  \end{displaymath}
  The middle square in the
  above diagram is a push-out because $- \diamond D(\alpha_2)$ and colimits
  preserve push-outs.
\end{proof}

\begin{lem}
\label{cralem2a}
  Let $(\C_0,\diamond,u)$ be a cocomplete symmetric monoidal category and
  suppose that for every $c \in \C_0$, the functor $c \diamond -$ preserves
  push-outs. 
  Let
  $f: X \to Y$ be a map of finite $G$-sets and let 
  $E \to F$ be a map in the category $(\Z\C)(X)$ of functors from the
  translation category of $X$ to  $\Z\C_0$.
  Consider the functor
  $f_\diamond : (\Z\C)(X) \to (\Z\C)(Y)$ with $f_\diamond(C)(y)(i) =
  (\underset{x \in f^{-1}(y)} \Diamond C(x))(i)$. 
  If the square
  \begin{displaymath}
  \begin{CD}
    E(x)(i) @>>> E(x)(i+1) \\
    @VVV @VVV \\ 
    F(x)(i) @>>> F(x)(i+1) 
  \end{CD}
  \end{displaymath}
  is a push-out in $\C$ for every $i \in \Z$ and $x \in X$, then the square
  \begin{displaymath}
  \begin{CD}
    (f_{\diamond}E)(y)(i)
    @>>> 
    (f_{\diamond}E)(y)(i+1) \\ 
    @VVV @VVV \\ 
    (f_{\diamond}F)(y)(i) @>>> (f_{\diamond}F)(y)(i+1).
  \end{CD}
  \end{displaymath}
  is a push-out in $\C$ for every $i \in \Z$ and $y \in Y$.
\end{lem}
\begin{proof}
  We have to check that for every $i \in \Z$ and $y \in Y$ the diagram
  \begin{displaymath}
  \begin{CD}
    (\underset {x \in f^{-1}(y)} \Diamond E(x))(i) 
    @>>>
    (\underset {x \in f^{-1}(y)} \Diamond E(x))(i+1) \\
    @VVV @VVV \\
    (\underset {x \in f^{-1}(y)} \Diamond F(x))(i) 
    @>>>
    (\underset {x \in f^{-1}(y)} \Diamond F(x))(i+1) 
  \end{CD}
  \end{displaymath}
  is a push-out diagram in $\C$. This follows from Lemma \ref{cralem2b}.
\end{proof}
The following
lemma can be used to identify filtration quotients of the form
$(D_1 \diamond D_2 \diamond \dots \diamond D_n)(k)/ (D_1 \diamond D_2
\diamond \dots \diamond D_n)(k-1)$ for filtered objects
$D_1,\dots,D_n$.
We leave its proof to the reader.    
\begin{lem}
\label{pufilt}
  Let $U$ be a finite set and consider, for $k \in \Z$, the sets
  $V_k \subseteq V_{\le k} \subseteq \Z^U$ 
  consisting of maps $\alpha : U \to \Z$ satisfying that 
  $\sum_{u \in U} \alpha(u) = k$ and $\sum_{u\in U} \alpha(u) \le k$
  respectively. 
  Let us consider these sets as partially ordered sets with $\beta \le \alpha$
  if and only if $\beta(u) \le \alpha(u)$ for every $u \in U$.
  Given 
  $\alpha \in \Z^U$ we let $V_{< \alpha}$ denote the partially ordered set
  consisting of those  
  $\beta \in \Z^U \setminus \{\alpha\}$
  satisfying that $\beta \le \alpha$.
  For every cocomplete category $\C$ and every
  functor $T: (\Z^U, \le) \to \C$ we have a push-out diagram of
  the form: 
  \begin{displaymath}
  \begin{CD}
    \coprod_{\alpha \in  V_k} \colim{\beta \in V_{<\alpha}} T(\beta) @>>> \coprod_{\alpha
    \in V_k} T(\alpha) \\
    @VVV @VVV \\
    \colim{\beta \in V_{\le k -1}} T(\beta) @>>>
    \colim{\alpha \in V_{\le k}} T(\alpha).
  \end{CD}
  \end{displaymath}
\end{lem}

\newpage
\noindent
\address{Morten Brun \\ Universit\"at Osnabr\"uck \\  Fachbereich Mathematik/Informatik \\
Albrechtstr. 28 \\
49069 Osnabr\"uck \\
Germany}

\vspace{0.5cm}
\noindent
   \email{brun@mathematik.uni-osnabrueck.de}
\end{document}